%
%
%
%
\hsize=5in
\baselineskip=12pt
\vsize=20.2cm
\parindent=10pt
\pretolerance=40
\predisplaypenalty=0
\displaywidowpenalty=0
\finalhyphendemerits=0
\hfuzz=2pt
\frenchspacing
\footline={\ifnum\pageno=1\else\hfil\tenrm\number\pageno\hfil\fi}
%
%
\input amssym.def
\def\titlefonts{\baselineskip=1.44\baselineskip
	\font\titlef=cmbx12
	\titlef
	}
\font\ninerm=cmr9
\font\ninebf=cmbx9
\font\ninei=cmmi9
\skewchar\ninei='177
\font\ninesy=cmsy9
\skewchar\ninesy='60
\font\nineit=cmti9
\def\reffonts{\baselineskip=0.9\baselineskip
	\textfont0=\ninerm
	\def\rm{\fam0\ninerm}%
	\textfont1=\ninei
	\textfont2=\ninesy
	\textfont\bffam=\ninebf
	\def\bf{\fam\bffam\ninebf}%
	\def\it{\nineit}%
	}
%
%
\def\frontmatter{\vbox{}\vskip1cm\bgroup
	\leftskip=0pt plus1fil\rightskip=0pt plus1fil
	\parindent=0pt
	\parfillskip=0pt
	\pretolerance=10000
	}
\def\endfrontmatter{\egroup\bigskip}
\def\title#1{{\titlefonts#1\par}}
\def\author#1{\bigskip#1\par}
\def\section#1\par{\ifdim\lastskip<\bigskipamount\removelastskip\fi
	\penalty-250\bigskip
	\vbox{\leftskip=0pt plus1fil\rightskip=0pt plus1fil
	\parindent=0pt
	\parfillskip=0pt
  \pretolerance=10000{\bf#1}}\nobreak\medskip
	}
\def\emph#1{{\it #1}\/}
\def\proclaim#1. {\medbreak\bgroup{\noindent\bf#1.}\ \it}
\def\endproclaim{\egroup
	\ifdim\lastskip<\medskipamount\removelastskip\medskip\fi}
\newdimen\itemsize
\def\setitemsize#1 {{\setbox0\hbox{#1\ }
	\global\itemsize=\wd0}}
\def\item#1 #2\par{\ifdim\lastskip<\smallskipamount\removelastskip\smallskip\fi
	{\leftskip=\itemsize
	\noindent\hskip-\leftskip
	\hbox to\leftskip{\hfil\rm#1\ }#2\par}\smallskip}
\def\Proof#1. {\ifdim\lastskip<\medskipamount\removelastskip\medskip\fi
	{\noindent\it Proof\if\space#1\space\else\ \fi#1.}\ }
\def\endproof{\hfill\hbox{}\quad\hbox{}\hfill\llap{$\square$}\medskip}
\def\Remark. {\ifdim\lastskip<\medskipamount\removelastskip\medskip\fi
        {\noindent\bf Remark. }}
\def\endremark{\medskip}
%
%
\newcount\citation
\newtoks\citetoks
\def\citedef#1\endcitedef{\citetoks={#1\endcitedef}}
\def\endcitedef#1\endcitedef{}
\def\citenum#1{\citation=0\def\curcite{#1}%
	\expandafter\checkendcite\the\citetoks}
\def\checkendcite#1{\ifx\endcitedef#1?\else
	\expandafter\lookcite\expandafter#1\fi}
\def\lookcite#1 {\advance\citation by1\def\auxcite{#1}%
	\ifx\auxcite\curcite\the\citation\expandafter\endcitedef\else
	\expandafter\checkendcite\fi}
\def\cite#1{\makecite#1,\cite}
\def\makecite#1,#2{[\citenum{#1}\ifx\cite#2]\else\expandafter\clearcite\expandafter#2\fi}
\def\clearcite#1,\cite{, #1]}
%
%
\def\references{\section References\par
	\bgroup
	\parindent=0pt
	\reffonts
	\rm
	\frenchspacing
	\setbox0\hbox{99. }\leftskip=\wd0
	}
\def\endreferences{\egroup}
\newtoks\authtoks
\newif\iffirstauth
\def\checkendauth#1{\ifx\auth#1%
    \iffirstauth\the\authtoks
    \else{} and \the\authtoks\fi,%
  \else\iffirstauth\the\authtoks\firstauthfalse
    \else, \the\authtoks\fi
    \expandafter\nextauth\expandafter#1\fi
	}
\def\nextauth#1,#2;{\authtoks={#1 #2}\checkendauth}
\def\auth#1{\nextauth#1;\auth}
\newif\ifinbook
\newif\ifbookref
\def\nextref#1 {\par\hskip-\leftskip
	\hbox to\leftskip{\hfil\citenum{#1}.\ }%
	\initnextref}
\def\initnextref{\bookreffalse\inbookfalse\firstauthtrue\ignorespaces}
\def\paper#1{{\it#1},}
\def\InBook#1{\inbooktrue in ``#1",}
\def\book#1{\bookreftrue{\it#1},}
\def\journal#1{#1\ifinbook,\fi}
\def\Vol#1{{\bf#1}}
\def\publisher#1{#1,}
\def\Year#1{\ifbookref #1.\else\ifinbook #1,\else(#1)\fi\fi}
\def\Pages#1{\makepages#1.}
\long\def\makepages#1-#2.#3{\ifinbook pp. \fi#1--#2\ifx\par#3.\fi#3}
\def\inRus{{ \rm(in Russian)}}
\def\etransl#1{English translation in \journal{#1}}
%
%
\newsymbol\square 1003
\newsymbol\smallsetminus 2272
\newsymbol\varnothing 203F
\newsymbol\bbk 207C
\newsymbol\barwedge 125A
\let\setm\smallsetminus
\let\wdg\barwedge

\def\defop#1#2{\def#1{\mathop{\rm #2}\nolimits}}
\defop\chr{char}
\defop\deg{deg}
\def\GL{GL}
\defop\HeckeSym{HeckeSym}
\defop\Id{Id}
\defop\Im{Im}
\defop\Ker{Ker}
\defop\rk{rank}
\def\mod#1{\mskip12mu(\mathop{\rm mod}#1)}

\def\Ypr#1#2{Y'\vphantom{Y}_{\!#1}^{\mskip1mu #2}}
\let\Rar\Rightarrow
\let\Lrar\Leftrightarrow
\let\ot\otimes
\let\sbs\subset
\let\al\alpha
\let\be\beta
\let\ga\gamma
\let\de\delta
\let\ep\varepsilon
\let\la\lambda
\let\om\omega
\def\omt{\widetilde\om}
\let\ph\varphi
\let\si\sigma
\let\ze\zeta
\let\De\Delta
\let\Up\Upsilon
\def\bbP{{\Bbb P}}
\def\bbS{{\Bbb S}}
\def\bbT{{\Bbb T}}
\def\bbZ{{\Bbb Z}}

\citedef
And17
Ar-TV90
Ar-TV91
Ew-O94
Gur90
Hua21
Hua22
Kl-Sch
Lyu86
Mo05
Shi23
Sk23a
Sk23b
Zh96
\endcitedef

\frontmatter
\title{Hecke symmetries associated with\break 
twisted polynomial algebras in 3 indeterminates}
\author{Nikita Shishmarov and Serge Skryabin}
\endfrontmatter

\medskip
{
\leftskip=10mm\rightskip=10mm\noindent
\reffonts
{\bf Abstract.}
\rm
We consider Hecke symmetries on a 3-dimensional vector space with the 
associated $R$-symmetric algebra isomorphic to the polynomial algebra 
$\bbk[x_1,x_2,x_3]$ twisted by an automorphism. The main result states that 
any such a Hecke symmetry is itself a twist of a Hecke symmetry with the 
associated $R$-symmetric algebra isomorphic to $\bbk[x_1,x_2,x_3]$. This 
allows us to describe equivalence classes of such Hecke symmetries.\par
}

\section
Introduction

With each Hecke symmetry $R$ on a vector space $V$ one associates the 
$R$-symmetric algebra $\bbS(V,R)$ which is viewed as a noncommutative analog 
of the ordinary symmetric algebra $\bbS(V)$ \cite{Gur90}. In a preceding 
article \cite{Sk23b} the second author gave a description of all Hecke 
symmetries on a 3-dimensional vector space $V$ such that $\bbS(V,R)=\bbS(V)$. 
The new paper extends that work by dealing with the cases where $\bbS(V,R)$ is 
isomorphic to the algebra $\bbS(V)$ twisted by an automorphism. The main 
result stated in Theorem 1.1 says that any such a Hecke symmetry is itself a 
twist of some Hecke symmetry with the respective $R$-symmetric algebra equal 
to $\bbS(V)$, and thus all such Hecke symmetries can be easily determined. We 
work over an arbitrary field $\bbk$, although one of arguments require 
$\chr\bbk\ne2$, and this is the only restriction in Theorem 1.1.

The question we pursue is opposite in a sense to the common direction in the 
study of Nichols algebras where one attempts to understand the structure of 
algebras associated with some classes of braidings (see a survey by 
Andruskiewitsch \cite{And17}). For the Hecke symmetries corresponding to 
quantum analogs of $\GL_3$ we know already that the respective algebras 
$\bbS(V,R)$ are Artin-Schelter regular of global dimension 3, and for each 
algebra in this class we want to find all braidings of Hecke type which yield 
the chosen algebra.

From the viewpoint of the classification of quantum $\GL_3$ groups splitting 
up the set of all relevant Hecke symmetries according to the algebra $\bbS(V,R)$ 
appears quite natural as this algebra is a primary invariant of a Hecke 
symmetry. This association of Hecke symmetries with algebras is not immediately 
retrievable from a paper of Ewen and Ogievetsky \cite{Ew-O94} where the 
classification is done up to twisting, and so some information gets lost. The 
approach discussed in \cite{Ew-O94} involves solving large systems of 
algebraic equations which would be difficult to analyze by hand, taking into 
account quite a large number of cases considered separately.

Our proof of Theorem 1.1 is accomplished in a fully invariant manner for an 
arbitrary twisting operator $\ze\in\GL(V)$. Generalizing what has been done in 
\cite{Sk23b} we associate with a Hecke symmetry $R$ a collection of linear 
forms $\ell_{xy}\in V^*$ indexed by pairs of vectors $x,y\in V$, and then 
reformulate conditions imposed on $R$ in terms of certain identities in the 
exterior algebra of the dual space $V^*$. However, the equations are now more 
complicated than those in \cite{Sk23b}, and we have to struggle to exclude the 
possibility of some other solutions. We apply geometric arguments to analyze 
identities satisfied by several polynomial maps arising in this context. The 
whole proof occupies sections 2 and 3 of the paper.

Explicit determination of Hecke symmetries with the prescribed $R$-symmetric 
algebra $\bbS(V,R)$ does depend on the twisting operator $\ze$. However, once 
it has been established that they are all twists, it requires only basic 
methods of linear algebra to select in the set of known Hecke symmetries with 
associated algebra $\bbS(V)$ those which commute with the linear operator 
$\ze\ot\ze$. We discuss parametrization of such Hecke symmetries in section 4 
of the paper.

Among all Artin-Schelter regular graded algebras of global dimension 3 the 
twisted polynomial algebras are characterized by the property that the scheme 
parametrizing their point modules is the whole projective plane $\bbP^2$ (see 
\cite{Ar-TV90}). Thus the present paper gives a description of all Hecke 
symmetries associated with algebras of this type. Unfortunately, the remaining 
Artin-Schelter regular algebras for which the point scheme is a cubic divisor 
in $\bbP^2$ defy unified treatment, and one has to resort to the case-by-case 
analysis, as outlined in \cite{Ew-O94}. For some cases details have been 
supplied in \cite{Shi23}, \cite{Sk23a}.

All our considerations are limited to the case $\,\dim V=3$. The quantum 
Yang-Baxter equation is notoriously hard to solve in higher dimensions. One of 
natural questions is whether the conclusion of Theorem 1.1 remains true when 
$\,\dim V>3$.

The twisted Hecke symmetries in our paper are related to some easily 
constructed 2-cocycles on the FRT bialgebras. Finding all 2-cocycles on the 
coordinate algebra of the group $\GL_n$ amounts to the determination of all 
quantum analogs of $\GL_n$. In its full generality this problem is 
intractable. In our modest contribution to the problem twisting is obtained by 
a simple conjugation in the group $\GL(V\ot V)$, and we do not need to ever 
invoke the FRT bialgebras.

\section
1. Twisting of graded algebras and Hecke symmetries

By the term ``algebra'' we understand a unital associative algebra over a 
field. Given a $\bbZ$-graded algebra $A=\bigoplus A_k$ and its automorphism 
$\si$ preserving the grading, the \emph{twist} of $A$ by $\si$ is the algebra 
$A_\si$ which has the same elements as the algebra $A$, but a different 
multiplication $\cdot_\si$ defined by the rule
$$
a\cdot_\si b=a\,\si^n(b)\qquad\hbox{for $\,a\in A_n$, $\,b\in A$}.
$$
(see \cite{Ar-TV91, section 8}).

If $A$ is generated by its component $V=A_1$, then every automorphism of $A$ 
is determined by its action on elements of degree 1. Let $A\cong\bbT(V)/I$ 
where $I=\bigoplus I_k$ is a graded ideal of the tensor algebra 
$\bbT(V)=\bigoplus_{k\ge0}V^{\ot k}$. An invertible linear map $\tau:V\to V$ 
extends to an automorphism of $A$ if and only if $\tau^{\ot k}(I_k)=I_k$ for 
all $k$. Let us denote by $A_\tau$ the twist of $A$ by the automorphism 
whose restriction to $A_1$ is $\tau$. This algebra is also generated by the 
degree 1 component. In fact, $A_\tau\cong\bbT(V)/I'$ where $I'$ is the graded 
ideal of $\bbT(V)$ with homogeneous components
$$
I'_k=(\Id_V\ot\,\tau^{-1}\ot\cdots\ot\tau^{-(k-1)})(I_k),\qquad k\ge2.
$$
In particular, $A_\tau$ has quadratic defining relations whenever so does $A$, 
and in this case the ideal $I'$ is generated by
$$
I'_2=(\Id_V\ot\,\tau^{-1})(I_2)=(\tau\ot\Id_V)(I_2).
$$

Denote by $\bbS(V)$ the symmetric algebra of a vector space $V$. It is the 
factor algebra of $\bbT(V)$ by the ideal generated by $\{x\ot y-y\ot x\mid 
x,y\in V\}$. Every invertible linear operator $\tau:V\to V$ extends to an 
automorphism of $\bbS(V)$, and we can form the respective twisted algebra 
$\bbS(V)_\tau$.

There are more general twistings of graded algebras introduced by Zhang 
\cite{Zh96}. However, in the case of $\bbS(V)$ every generalized twist in the 
sense of Zhang is isomorphic as a graded algebra to the twist of $\bbS(V)$ by 
an automorphism (see \cite{Zh96, Prop. 5.13}). Twists by automorphisms can be 
realized as a very special case of cocycle twists of comodule algebras 
discussed, e.g., in \cite{Kl-Sch, 10.3.2} and \cite{Mo05}. With respect to 
2-cocycles on Manin's universally coacting Hopf algebra every Artin-Schelter 
regular graded algebra with the Hilbert series $1/(1-t)^n$ is a twist of the 
polynomial algebra in $n$ indeterminates (see \cite{Hua22, Th. 5.1.1}).

Let $V$ be a finite-dimensional vector space over a field $\bbk$. A 
\emph{Hecke symmetry} with parameter $0\ne q\in\bbk$ on this vector space is 
a linear operator $R:V\ot V\to V\ot V$ satisfying the braid equation
$$ 
(R\ot\Id_V)(\Id_V\ot\,R)(R\ot\Id_V)=(\Id_V\ot\,R)(R\ot\Id_V)(\Id_V\ot\,R)\eqno(1.1)
$$
and the quadratic Hecke relation
$$
(R-q\cdot\Id_{V\ot V})(R+\Id_{V\ot V})=0.\eqno(1.2)
$$
This notion is due to Gurevich \cite{Gur90}, while involutive symmetries 
corresponding to the parameter value $q=1$ were studied earlier by 
Lyubashenko \cite{Lyu86}.

The \emph{$R$-symmetric algebra} $\bbS(V,R)$ associated with $R$ is the factor 
algebra of $\bbT(V)$ by the graded ideal generated by the subspace
$$ 
\Im\,(R-q\cdot\Id_{V\ot V})\sbs V^{\ot 2}.
$$

If $\ze:V\to V$ is an invertible linear operator such that $\ze\ot\ze$ commutes 
with $R$, then the linear operator
$$
R_\ze=(\ze\ot\Id_V)\circ R\circ(\ze^{-1}\ot\Id_V)
=(\Id_V\ot\,\ze^{-1})\circ R\circ(\Id_V\ot\,\ze)\eqno(1.3)
$$
is also a Hecke symmetry called the \emph{twist} of $R$ by $\ze$. It satisfies 
the Hecke relation with the same parameter $q$, and the braid equation for 
$R_\ze$ follows from the equalities
$$
\openup1\jot
\eqalign{ 
(R_\ze\ot\Id_V)&(\Id_V\ot\,R_\ze)(R_\ze\ot\Id_V)\cr
&=(\ze\ot\Id_V\ot\,\ze^{-1})(R\ot\Id_V)(\Id_V\ot\,R)(R\ot\Id_V)(\ze^{-1}\ot\Id_V\ot\,\ze),\cr
(\Id_V\ot\,R_\ze)&(R_\ze\ot\Id_V)(\Id_V\ot\,R_\ze)\cr
&=(\ze\ot\Id_V\ot\,\ze^{-1})(\Id_V\ot\,R)(R\ot\Id_V)(\Id_V\ot\,R)(\ze^{-1}\ot\Id_V\ot\,\ze).
}
$$
Such twists were used by Ewen and Ogievetsky \cite{Ew-O94}. The condition that 
$\ze\ot\ze$ commutes with $R$ actually allows one to define a certain 2-cocycle 
on the FRT bialgebra associated with $R$. In coordinate form this 
condition is written out in \cite{Hua21, Lemma 3.1.4}. The FRT bialgebra 
associated with $R_\ze$ is a cocycle twist of the initial bialgebra. However, 
we will never need this interpretation.

The ideal of the tensor algebra $\bbT(V)$ defining the $R_\ze$-symmetric 
algebra $\bbS(V,R_\ze)$ is generated by the space
$$
\Im(R_\ze-q\cdot\Id_{V\ot V})=(\ze\ot\Id_V)\bigl(\Im(R-q\cdot\Id_{V\ot V})\bigr).
$$
Hence $\,\,\bbS(V,R_\ze)\cong\bbS(V,R)_\ze$.

The initial Hecke symmetry $R$ is itself a twist of $R_\ze$. Indeed, 
$\ze\ot\ze$ commutes with $R_\ze$, and $\,R=(R_\ze)_{\ze^{-1}}$.\quad 
Our main result is

\proclaim
Theorem 1.1.
Let $V$ be a vector space of dimension $3$ over a field $\,\bbk$ of 
characteristic $\ne2$, and let $\ze\in\GL(V)$. If $R$ is a Hecke symmetry on 
$V$ such that the identity transformation of $\,V$ extends to an isomorphism 
of algebras $\,\bbS(V,R)\cong\bbS(V)_\ze,$ then $R$ commutes with $\ze\ot\ze$.

As a consequence, such a Hecke symmetry is a twist of a Hecke symmetry with 
the associated $R$-symmetric algebra equal to $\,\bbS(V)$.
\endproclaim

The final assertion of Theorem 1.1 is an obvious consequence of the first one. 
If $R$ commutes with $\ze\ot\ze$, then $R_{\ze^{-1}}$ is a Hecke symmetry with 
the associated algebra $\bbS(V,R_{\ze^{-1}})=\bbS(V)$, and we have 
$\,R=(R_{\ze^{-1}})_\ze$.

In section 2 we will associate with the Hecke symmetry $R$ a collection 
of linear forms $\ell_{xy}\in V^*$ indexed by pairs of vectors $x,y\in V$ and 
show that the braid equation for $R$ implies that these linear forms satisfy 
a certain identity in the second exterior power $\bigwedge^2V^*$ of the dual 
space $V^*$. Likewise the property that $R$ commutes with $\ze\ot\ze$ will be 
reformulated in terms of a certain identity in $\bigwedge^2V^*$.

The proof of Theorem 1.1 will be completed in section 3 by analyzing all the 
identities found in section 2. The most difficult part presented in 
Proposition 3.2 excludes the possibility of one case which could lead to Hecke 
symmetries violating the conclusion of Theorem 1.1. The assumption 
$\chr\bbk\ne2$ will be used in the proof of Theorem 1.1 only once in Lemma 3.10.

\medskip
The twisting transformation $R\mapsto R_\ze$ is conjugation in the group 
$\GL(V\ot V)$ by means of the linear operator $\ze\ot\Id_V$. If $R$ is any 
Hecke symmetry on $V$, then for any $\ph\in\GL(V)$ the linear operator
$$
R'=(\ph\ot\ph)\circ R\circ(\ph\ot\ph)^{-1}\eqno(1.4)
$$
is again a Hecke symmetry on $V$, and $\ph$ extends to an isomorphism of 
graded algebras $\,\bbS(V,R)\to\bbS(V,R')$. Conjugation by the operators 
$\ph\ot\ph$ defines an action of the group $\GL(V)$ on the set of all Hecke 
symmetries on $V$ which we call the \emph{conjugating action}. We say that two 
Hecke symmetries on $V$ are \emph{equivalent} if they lie in the same orbit 
with respect to the conjugating action of $\GL(V)$.

If $R$ and $R'$ are two Hecke symmetries such that $\,\bbS(V,R)=\bbS(V,R')$ in 
the sense that the two algebras have the same space of defining relations in 
$V^{\ot2}$, then equality (1.4) implies that $\ph$ extends to an automorphism 
of $\,\bbS(V,R)$. In other words, two Hecke symmetries with the same 
associated $R$-symmetric algebra $A$ are equivalent if and only if they lie in 
the same orbit with respect to the conjugating action of the subgroup of 
$\GL(V)$ consisting of all linear operators on $V$ which extend to an 
automorphism of $A$.

In section 4 we will use Theorem 1.1 to describe equivalence classes of Hecke 
symmetries with the associated $R$-symmetric algebra isomorphic to 
$\,\bbS(V)_\ze$.

\section
2. Setup for the proof of the main result

Let $V$ be a vector space of dimension 3 over a field $\bbk$. 
We fix a linear operator $\ze\in\GL(V)$. Let $A=\bbT(V)/I$ where $I$ is 
the graded ideal of $\bbT(V)$ generated by the set 
$\{\ze(x)\ot y-\ze(y)\ot x\mid x,y\in V\}$. So $A\cong\bbS(V)_\ze$.

Further on we will write elements of $\bbT(V)$ omitting the sign $\ot$, i.e.,
$xy=x\ot y$, $\,xyz=x\ot y\ot z$ for $x,y,z\in V$. Put
$$
\eqalign{
x\wdg y={}&\ze(x)\,y-\ze(y)\,x\in V^{\ot2}\qquad\hbox{and}\cr
\noalign{\smallskip}
x\wdg y\wdg z={}&
\ze^2(x)\,\ze(y)\,z+\ze^2(y)\,\ze(z)\,x+\ze^2(z)\,\ze(x)\,y\cr
&-\ze^2(x)\,\ze(z)\,y-\ze^2(y)\,\ze(x)\,z-\ze^2(z)\,\ze(y)\,x\in V^{\ot3}.
}
$$
Consider the graded subspace $\Up=\bigoplus\Up^{(k)}$ of $\bbT(V)$ where
$$
\Up^{(0)}=\bbk,\qquad\Up^{(1)}=V,\qquad\Up^{(2)}=I_2,\qquad
\Up^{(3)}=I_2V\cap VI_2
$$
and $\Up^{(k)}=0$ for $k>3$. Its components $\Up^{(2)}$ and $\Up^{(3)}$ are 
spanned by the tensors, respectively, $x\wdg y$ and $x\wdg y\wdg z$. There are 
linear isomorphisms $\bigwedge^kV\cong\Up^{(k)}$ given by the identity maps 
for $k=0$ and $k=1$, and such that
$$
x\wedge y\mapsto x\wdg y,\qquad x\wedge y\wedge z\mapsto x\wdg y\wdg z.
$$
for $k=2$ and $k=3$. By means of these isomorphisms we obtain an algebra 
structure on $\Up$, and $\wdg$ can be understood as the respective 
multiplication. Thus $(\Up,\wdg)$ is a graded Frobenius algebra isomorphic to 
the exterior algebra of $V$. The sole purpose of this multiplication on $\Up$ 
is that it gives a nondegenerate bilinear pairing between $V$ and $I_2$ which 
will be used in formulas (2.8).

Note that
$$
\ze^{\ot3}(x\wdg y\wdg z)=(\det\ze)\,x\wdg y\wdg z\eqno(2.1)
$$
since this tensor is the image of 
$\,\ze(x)\wedge\ze(y)\wedge\ze(z)\in\bigwedge^3V$.

Suppose that $R$ is a Hecke symmetry on $V$ such that $\bbS(V,R)=A$. In other 
words, $R$ satisfies the hypothesis of Theorem 1.1. Put
$$
R'=(\ze^{-1}\ot\ze^{-1})\circ R\circ(\ze\ot\ze).\eqno(2.2)
$$
Then $R'=R$ if and only if $\ze\ot\ze$ commutes with $R$. In any case $R'$ is 
a Hecke symmetry with the same associated algebra $\bbS(V,R')=A$ since 
$(\ze\ot\ze)(I_2)=I_2$ (so $\ze$ extends to an automorphism of $A$).

Let $q$ be the parameter of the Hecke relation, and
$$
Y=q\cdot\Id_{V\ot V}-\,R\eqno(2.3)
$$
the \emph{$R$-skewsymmetrizer}. Then $\,\Im Y=\Up^{(2)}$ and $Y^2=(q+1)Y$. Hence
$$
Yw=(q+1)w\quad\hbox{for all $w\in\Up^{(2)}$}.\eqno(2.4)
$$
The braid equation for $R$ can be rewritten as the equation
$$
\eqalign{
&(\Id_V\ot\,Y)(Y\ot\Id_V)(\Id_V\ot\,Y)-q\cdot(\Id_V\ot\,Y)\cr
&\hskip30mm=(Y\ot\Id_V)(\Id_V\ot\,Y)(Y\ot\Id_V)-q\cdot(Y\ot\Id_V).
}\eqno(2.5)
$$
Here the linear operators $\Id\ot Y$ and $Y\ot\Id$ acting on $V^{\ot3}$ have 
images, respectively, $V\ot\Up^{(2)}$ and $\Up^{(2)}\ot V$. Therefore the two 
equal operators in (2.5) have images in the 1-dimensional subspace
$$
\Up^{(3)}=(V\ot\Up^{(2)})\cap(\Up^{(2)}\ot V)\sbs V^{\ot3}.
$$
So it follows that
$$
(\Id_V\ot\,Y)(Y\ot\Id_V)w-qw\in\Up^{(3)}\quad\hbox{for all 
$w\in V\ot\Up^{(2)}$}.\eqno(2.6)
$$

Fix a nonzero alternating trilinear form $\om:V\times V\times V\to\bbk$. 
There is a linear bijection $\omt:\Up^{(3)}\to\bbk$ such that
$$
\omt(x\wdg y\wdg z)=\om(x,y,z)\quad\hbox{for $x,y,z\in V$}.\eqno(2.7)
$$
Define linear forms $\,\ell_{xy}\,,\,\ell'_{xy}\in V^*\,$ by the rule
$$
\openup1\jot
\eqalign{
\ell_{xy}(z)&=\omt\bigl(x\wdg Y\bigl(\ze(y)z\bigr)\bigr),\cr
\ell'_{xy}(z)&=\omt\bigl(x\wdg Y'\bigl(\ze(y)z\bigr)\bigr),\qquad 
x,y,z\in V,
}\eqno(2.8)
$$
where $Y'=q\cdot\Id-R'$ is the $R'$-skewsymmetrizer. In this way we obtain two 
collections of linear forms on $V$ indexed by pairs of vectors $x,y\in V$. 
Both $\ell_{xy}$ and $\ell'_{xy}$ depend linearly on $x$ and on $y$. Note that
$$
Y'=(\ze^{-1}\ot\ze^{-1})\circ Y\circ(\ze\ot\ze)\eqno(2.9)
$$
and therefore
$$
\ell'_{xy}(z)=(\det\ze)^{-1}\,\ell_{\ze(x)\ze(y)}\bigl(\ze(z)\bigr).\eqno(2.10)
$$

Containments (2.6) are not fully equivalent to the braid equation for $R$, 
but they are sufficient for our purposes. What is important, (2.6) can be 
interpreted in terms of the functions $\ell_{xy}$ and $\ell'_{xy}$. We will 
not need the most general identity equivalent to (2.6) and will use only its 
special case stated below:

\proclaim
Lemma 2.1.
Containments $(2.6)$ imply the following identity in $\,\bigwedge^2V^*:$
$$
\ell_{xy}\wedge\ell'_{xy}=\ell_{xx}\wedge\ell'_{yy}\,,\qquad x,y\in V.\eqno(2.11)
$$
\endproclaim

\Proof.
For $f,g\in V^*$ we identify $f\wedge g\in\bigwedge^2V^*$ with an alternating 
bilinear form on $V$ setting
$$
(f\wedge g)(u,v)=f(u)\,g(v)-f(v)\,g(u),\qquad u,v\in V.\eqno(2.12)
$$

Let $e_1,e_2,e_3$ be any linear basis for the vector space $V$. In $V^{\ot2}$ 
and $V^{\ot3}$ we take the bases, respectively, $\{\ze(e_i)\,e_j\}$ and 
$\{\ze^2(e_i)\,\ze(e_j)\,e_k\}$. They are better adapted to $\Up$ than 
the natural bases $\{e_ie_j\}$ and $\{e_ie_je_k\}$. Let 
$Y_{ij}^{st},\,\Ypr{ij}{st}\in\bbk$ be the coefficients in the expressions
$$
Y\bigl(\ze(e_i)\,e_j\bigr)=\sum Y_{ij}^{st}\,\ze(e_s)\,e_t,\qquad 
Y'\bigl(\ze(e_i)\,e_j\bigr)=\sum\Ypr{ij}{st}\,\ze(e_s)\,e_t.
$$
Here and later it is assumed that the summation is over the indices repeated 
as subscripts and superscripts. Since $Y$ and $Y'$ have images in $\Up^{(2)}$, 
we have
$$
Y_{ij}^{st}=-Y_{ij}^{ts},\quad Y_{ij}^{ss}=0,\quad 
\Ypr{ij}{st}=-\Ypr{ij}{ts},\quad\Ypr{ij}{ss}=0.
$$
Now
$$
Y\bigl(\ze^2(e_i)\,\ze(e_j)\bigr)=(\ze\ot\ze)\bigl(Y'\bigl(\ze(e_i)\,e_j\bigr)\bigr)
=\sum\Ypr{ij}{rl}\,\ze^2(e_r)\,\ze(e_l),
$$
and for $\,w=\ze^2(e_i)\bigl(\ze(e_j)e_k-\ze(e_k)e_j\bigr)\in 
V\ot\Up^{(2)}\,$ we find
$$ 
(\Id\ot\,Y)(Y\ot\Id)(w)=\sum\bigl(\Ypr{ij}{rl}\,Y_{lk}^{st}-\Ypr{ik}{rl}\,Y_{lj}^{st}\bigr)\,
\ze^2(e_r)\ze(e_s)e_t.
$$
If $s=r$, then the basis element $\ze^2(e_r)\ze(e_s)e_t$ has zero coefficient 
in the elements of $\Up^{(3)}$, and if $i\ne r$, then $\ze^2(e_r)\ze(e_s)e_t$ 
has zero coefficient in $w$. Therefore (2.6) implies that
$$
\sum\,\bigl(\Ypr{ij}{rl}\,Y_{lk}^{rt}-\Ypr{ik}{rl}\,Y_{lj}^{rt}\bigr)=0
\quad\hbox{if $i\ne r$}\eqno(2.13)
$$
where the summation is over $l$. Now note that
$$
\openup1\jot
\displaylines{
Y_{ij}^{23}=\ell_{e_1e_i}(e_j)\mskip2mu\al,\qquad
Y_{ij}^{31}=\ell_{e_2e_i}(e_j)\mskip2mu\al,\qquad
Y_{ij}^{12}=\ell_{e_3e_i}(e_j)\mskip2mu\al,\cr
\Ypr{ij}{23}=\ell'_{e_1e_i}(e_j)\mskip2mu\al,\qquad
\Ypr{ij}{31}=\ell'_{e_2e_i}(e_j)\mskip2mu\al,\qquad
\Ypr{ij}{12}=\ell'_{e_3e_i}(e_j)\mskip2mu\al
}
$$
where $\al=\om(e_1,e_2,e_3)^{-1}$. These formulas follow straightforwardly from 
(2.8). For example,
$$
\ell_{e_1e_i}(e_j)=\omt\bigl(e_1\wdg Y\bigl(\ze(e_i)e_j\bigr)\bigr)
=Y_{ij}^{23}\,\om(e_1,e_2,e_3)=Y_{ij}^{23}\,\al^{-1}
$$
since $\,\,Y\bigl(\ze(e_i)e_j\bigr)
=Y_{ij}^{23}\,e_2\wdg e_3+Y_{ij}^{31}\,e_3\wdg e_1+Y_{ij}^{12}\,e_1\wdg e_2\,$.

If $r=t$ then each term in the left hand side of (2.13) vanishes. Suppose 
$r\ne t$, and let $p\in\{1,2,3\}$ be the remaining element $\ne r,t$. 
Making use of (2.12) we get
$$
\Ypr{ij}{rl}\,Y_{lk}^{rt}-\Ypr{ik}{rl}\,Y_{lj}^{rt}=\cases{
\hphantom{\;(\ell_{e_ie_i}\wedge\ell'_{e_te_i})(e_j,e_k)}0 & for $l=r$\cr
(\ell_{e_pe_t}\wedge\ell'_{e_pe_i})(e_k,e_j)\,\be^2 & for $l=t$\cr
(\ell_{e_pe_p}\wedge\ell'_{e_te_i})(e_j,e_k)\,\be^2 & for $l=p$
}
$$
where $\be=\om(e_p,e_r,e_t)^{-1}$. Take $i=t$. So $i\ne r$, and (2.13) can be 
rewritten as
$$
(\ell_{e_pe_t}\wedge\ell'_{e_pe_t}-\ell_{e_pe_p}\wedge\ell'_{e_te_t})(e_k,e_j)=0.
$$
Since $e_k$ and $e_j$ are two arbitrary basis vectors, the above equality is 
exactly (2.11) with $x=e_p$ and $y=e_t$.

All considerations are valid with respect to any basis of $V$. If $x,y\in V$ 
are two linearly independent vectors, we can include them in some basis of 
$V$, and then (2.11) for these vectors $x,y$ will follow by making use of the 
chosen basis. On the other hand, if both $x$ and $y$ are scalar multiples 
of some vector, then equality (2.11) is automatically true.
\endproof

\Remark.
By analyzing equations more fully one can see that for the containments $(2.6)$ 
to hold it is necessary and sufficient that the equality
$$
(\ell_{xy}\wedge\ell'_{xz}-\ell_{xx}\wedge\ell'_{yz})(u,v)=q\,\om(x,y,z)\,\om(x,u,v)
$$
be true for all $x,y,z,u,v\in V$. This generalizes the case $\ze=\Id_V$ 
considered in \cite{Sk23b, Lemma 2.2}. If $\ze=\Id_V$, then 
$\ell'_{xy}=\ell_{xy}$, and (2.11) reduces to a much simpler identity 
$\ell_{xx}\wedge\ell_{yy}=0$ which was used in \cite{Sk23b} to determine all 
Hecke symmetries with the associated $R$-symmetric algebra 
$\,\bbS(V,R)=\bbS(V)$.
\endremark

Equality (2.4) too can be interpreted in terms of the functions $\ell_{xy}$. 
Indeed, it means that $\,\omt(x\wdg Yw)=(q+1)\,\omt(x\wdg w)\,$ for all 
$x\in V$ and $w\in\Up^{(2)}$. Taking $w=y\wdg z=\ze(y)\,z-\ze(z)\,y$, we 
rewrite this as
$$
\ell_{xy}(z)-\ell_{xz}(y)=(q+1)\,\om(x,y,z)\,,\qquad x,y,z\in V.
$$
Since $R'$ is a Hecke symmetry with the same associated algebra 
$\bbS(V,R')=\bbS(V,R)$, there is a similar identity with the functions 
$\ell_{xy}$ replaced by $\ell'_{xy}$. Again, we will use only a special 
case of these identities obtained by taking $z=x$:

\proclaim
Lemma 2.2.
The condition that the linear operators $Y$ and $Y'$ act on the subspace 
$\Up^{(2)}\sbs V^{\ot2}$ as $(q+1)$ times the identity operator implies that
$$
\ell_{xy}(x)=\ell_{xx}(y),\qquad\ell'_{xy}(x)=\ell'_{xx}(y)\eqno(2.14)
$$
for all $x,y\in V$.
\endproclaim

The braid equation implies other useful properties of the collection 
$\{\ell_{xy}\}$. They are immediate consequences of the following fact:

\proclaim
Lemma 2.3.
If some vector $a\in V$ has the property that either $Y(ax)=0$ for all 
$x\in V$ or $Y(xa)=0$ for all $x\in V$, then $a=0$.
\endproclaim

\Proof.
The equality $Y(ax)=0$ means that $R(ax)=qax$. If this holds for all $x\in V$ 
then, applying the operators in the braid equation (1.1) to the tensor 
$w=axy\in V^{\ot3}$, we get
$$
q^2\,(\Id_V\ot\,R)w=q\,(\Id_V\ot\,R^2)w.
$$
If $a\ne0$, then it follows that $qR(xy)=R^2(xy)$ for all $x,y\in V$, whence 
$R=q\cdot\Id$ since $R$ is an invertible operator. Thus $Y=0$ by (2.3). 
However, this contradicts the assumption that $\,\Im Y=\Up^{(2)}\,$ is a space 
of dimension 3.

In the other case when $Y(xa)=0$ for all $x\in V$ we argue similarly, now 
working with the tensors of the form $\,w=xya$.
\endproof

\proclaim
Lemma 2.4.
The set $\,\{\ell_{xy}\mid x,y\in V\}\,$ spans the whole space $V^*$. If $a\in 
V$ is such that either $\ell_{ax}=0$ for all $x\in V$ or $\ell_{xa}=0$ for all 
$x\in V,$ then $a=0$.
\endproclaim

\Proof.
For the first assertion we have to check that the only vector at which all 
linear forms $\ell_{xy}$ vanish is the zero vector. But $\ell_{xy}(a)=0$ means 
that $x\wdg Y\bigl(\ze(y)a\bigr)=0$, and if this equality holds for all 
$x,y\in V$, then $Y(va)=0$ for all $v\in V$. So Lemma 2.3 applies.

If $\ell_{xa}=0$ for all $x\in V$, then $Y\bigl(\ze(a)v\bigr)=0$ for all 
$v\in V$. If $\ell_{ax}=0$ for all $x\in V$, then $a\wdg t=0$ for all 
$t\in\Up^{(2)}$. The conclusion $a=0$ follows in each case.
\endproof

Our aim is to show that $R$ commutes with $\ze\ot\ze$, and we are going to 
reformulate this property in terms of the functions $\ell_{xy}$, $\ell'_{xy}$. 
This will be done in Lemma 2.6. First we mention one related result derived 
from an earlier work:

\proclaim
Lemma 2.5.
Any Hecke symmetry $R$ satisfying the hypothesis of Theorem $1.1$ commutes with 
the linear operator $(\ze\ot\ze)^3$.
\endproclaim

\Proof.
The tensors $t$ in the one-dimensional subspace $\Up^{(3)}\sbs V^{\ot3}$ 
satisfy the twisted cyclicity condition with the twisting operator 
$\psi=(\det\ze)\,\ze^{-3}$. This means that
$$
t=(\Id_V\ot\Id_V\ot\,\psi)\,s_{123}(t)
$$
where $s_{123}$ is the permutation operator on $V^{\ot3}$ defined by the rule 
$s_{123}(xyz)=yzx$ for $x,y,z\in V$. Since $t=x\wdg y\wdg z$ for suitable 
$x,y,z$, the displayed equality is readily verified with the aid of (2.1). By 
\cite{Sk23a, Theorem 3.8} $R$ commutes with $\psi\ot\psi$, and this gives the 
desired conclusion.
\endproof

Under certain conditions on the eigenvalues of $\ze$ the linear operators 
$\ze\ot\ze$ and $(\ze\ot\ze)^3$ have equal centralizers in the group 
$\GL(V\otimes V)$. However, in general the conclusion of Lemma 2.5 is weaker 
than what we need, and we do not rely on it. Lemma 2.5 will be used only once 
to provide special treatment of the case $q=-1$ in Lemma 2.6. For $q\ne-1$ 
Lemma 2.5 will not be used.

\setitemsize(iii)
\proclaim
Lemma 2.6.
The following conditions are equivalent:

\item(i)
$R$ commutes with $\ze\ot\ze,$

\item(ii)
$Y'=Y,$

\item(iii)
$Y'=cY$ for some $c\in\bbk,$

\item(iv)
there is $c\in\bbk$ such that $\,\ell'_{xy}=c\,\ell_{xy}\,$ for all $x,y\in V,$

\item(v)
$\ell_{xy}\wedge\ell'_{xy}=0$ for all $x,y\in V$.

\endproclaim

\Proof.
It is immediate from (2.2) that (i) holds if and only if $R'=R$, which in turn 
is equivalent to the equality $Y'=Y$ by the definition of these operators. 
We have to show that $c=1$ whenever (iii) holds.

So suppose that $Y'=cY$. Recall that both $Y$ and $Y'$ act on the subspace 
$\Up^{(2)}$ as $(q+1)$ times the identity operator. Hence $c(q+1)=q+1$. If 
$q\ne-1$, we do get $c=1$.

If $q=-1$, we have to argue differently. Note that containments (2.6) are 
valid also for $Y'$ since $R'$ is a Hecke symmetry with the same associated 
algebra $S(V,R')=A$. However, the assumption $Y'=cY$ implies that
$$ 
(\Id_V\ot\,Y')(Y'\ot\Id_V)w\equiv c^2qw\mod{\Up^{(3)}}\qquad\hbox{for all 
$w\in V\ot\Up^{(2)}$}.
$$
Hence $c^2q=q$, i.e., $c^2=1$. We now make use of Lemma 2.5. Since the linear 
operator $(\ze\ot\ze)^3$ commutes with $R$, it also commutes with $Y$. On the 
other hand, from (2.9) we deduce that
$$
Y\circ(\ze\ot\ze)=c\,(\ze\ot\ze)\circ Y.
$$
It follows that $c^3=1$. Together with $c^2=1$ this yields $c=1$.

Recalling (2.8), we see that (iv) is expanded as
$$
\omt\bigl(x\wdg Y'\bigl(\ze(y)z\bigr)\bigr)
=\omt\bigl(x\wdg c\,Y\bigl(\ze(y)z\bigr)\bigr)\quad\hbox{for all $x,y,z\in V$}.
$$
Since $\omt$ induces a nondegenerate bilinear pairing $V\times\Up^{(2)}\to\bbk$, 
this identity is equivalent to
$$
Y'\bigl(\ze(y)z\bigr)=c\,Y\bigl(\ze(y)z\bigr)\quad\hbox{for all $y,z\in V$},
$$
which amounts to (iii). Thus 
$\hbox{(i)}\Lrar\hbox{(ii)}\Lrar\hbox{(iii)}\Lrar\hbox{(iv)}$. It is also 
clear that (v) follows from (iv). The opposite implication is more 
complicated.

Suppose that (v) holds. If $x,y\in V$ are any two vectors such that 
$\ell_{xy}\ne0$, then (v) implies that $\ell'_{xy}$ is a scalar multiple of 
$\ell_{xy}$, and so $\ell'_{xy}=c(x,y)\,\ell_{xy}$ for some $c(x,y)\in\bbk$. 
We have to deal with the dependency of $c(x,y)$ on $x$ and $y$.

For each $x\in V$ put $K_x=\{y\in V\mid\ell_{xy}=0\}$. This is a vector 
subspace of $V$. By Lemma 2.4 we have $K_x\ne V$ unless $x=0$. Put
$$
X=\{x\in V\mid\dim K_x\le1\}.
$$

\proclaim
Claim 1. For each $x\in X$ there exists $c_\flat(x)\in\bbk$ such that 
$\,\ell'_{xy}=c_\flat(x)\,\ell_{xy}\,$ for all $y\in V$ simultaneously.
\endproclaim

Note that $K_x$ is the kernel of the linear map $V\to V^*$ given by the 
assignment $y\mapsto\ell_{xy}$. It follows that $\ell_{xy_1}$ and 
$\ell_{xy_2}$ are linearly independent linear forms whenever $y_1,y_2\in V$ 
are two vectors which are linearly independent modulo the one-dimensional 
subspace $K_x$. Then $\ell_{x(y_1+y_2)}=\ell_{xy_1}+\ell_{xy_2}\ne0$ and
$$
\ell'_{x(y_1+y_2)}=\ell'_{xy_1}+\ell'_{xy_2}=c(x,y_1)\,\ell_{xy_1}+c(x,y_2)\,\ell_{xy_2}.
$$
On the other hand, the left hand side is equal to 
$c(x,y_1+y_2)\,\ell_{x(y_1+y_2)}$. Hence
$$
c(x,y_1)=c(x,y_1+y_2)=c(x,y_2).
$$
It follows that there exists $c_\flat(x)\in\bbk$ such that $c(x,y)=c_\flat(x)$ for all 
$y\in V$ such that $y\notin K_x$. Now the set
$$
V_x=\{y\in V\mid\ell'_{xy}=c_\flat(x)\,\ell_{xy}\}
$$
is a vector subspace of $V$ which contains the set $V\setm K_x$. Since the 
latter spans the whole space $V$, we conclude that $V_x=V$.

\proclaim
Claim 2. The set $X$ is nonempty.
\endproclaim

Suppose that $X=\varnothing$. Then $\dim K_x>1$ for each $x\in V$. If $x\ne0$, 
then we must have $\dim K_x=2$ since $K_x\ne V$, and therefore the set 
$T_x=\{\ell_{xy}\mid y\in V\}$ is a one-dimensional subspace of $V^*$. Pick 
some nonzero vectors $x_1\in V$ and $a\in K_{x_1}$. The set 
$$
L=\{x\in V\mid\ell_{xa}=0\}
$$
is a vector subspace of $V$, and $L\ne V$ by Lemma 2.4. Hence $V\setm L$ spans 
the whole space $V$, and by Lemma 2.4 the linear span of the set 
$\{\ell_{xy}\mid x,y\in V,\ x\notin L\}$ must coincide with the whole $V^*$ 
since it contains the functions $\ell_{xy}$ for all $x,y\in V$. We can 
find $x_2\in V\setm L$ such that $T_{x_2}\ne T_{x_1}$.

Since $\ell_{x_2a}\ne0$, we have $a\notin K_{x_2}$. Hence $K_{x_1}\ne 
K_{x_2}$, while these two spaces both have dimension 2. Take some vectors 
$y_1\in K_{x_2}\setm K_{x_1}$ and $y_2\in K_{x_1}\setm K_{x_2}$. Then 
$\ell_{x_1y_1}$ and $\ell_{x_2y_2}$ are nonzero, and moreover these two 
functions are linearly independent since they lie in two distinct 
one-dimensional subspaces $T_{x_1}$ and $T_{x_2}$. On the other hand, 
$\ell_{x_1y_2}=\ell_{x_2y_1}=0$. Hence
$$
\ell_{(x_1+x_2)y_1}=\ell_{x_1y_1}\,,\qquad
\ell_{(x_1+x_2)y_2}=\ell_{x_2y_2}\,,
$$
and so $\,\dim T_{x_1+x_2}\ge2$, in contradiction with the assumption 
$\,\dim K_{x_1+x_2}>1$.

\proclaim
Claim 3. The linear span $\bbk X$ of the set $X$ coincides with $V$.
\endproclaim

Suppose that $\bbk X\ne V$. We have $\dim K_u=2$ for all $u\in V$ such that 
$u\notin\bbk X$. If $u\notin\bbk X$ and $x\in X$, then $u+x\notin\bbk X$, 
whence $\dim K_{u+x}=2$ as well. Furthermore, each element $a\in K_u\cap K_{u+x}$ 
is contained in $K_x$ since $\ell_{xa}=\ell_{(u+x)a}-\ell_{ua}=0$. This forces 
$K_u\ne K_{u+x}$ since $\dim K_x\le1$. Hence $K_u\cap K_{u+x}$ is a 
one-dimensional subspace of $V$, and we must have $K_x=K_u\cap K_{u+x}$.

It follows that $K_x\sbs K_u$ for each $x\in X$ and each $u\in V\setm\bbk X$. 
Moreover, if we put $K=\bigcap_{u\notin\bbk X}K_u$, then $K_x=K$ for each 
$x\in X$ and $K\sbs K_u$ for each $u\in V\setm\bbk X$. We have seen already 
that $K_x\ne0$ for $x\in X$. Hence $K\ne0$. If $a\in K$, then $\ell_{xa}=0$ 
for all $x\in X$, as well as for all $x\in V\setm\bbk X$. This equality holds 
then also for all $x\in\bbk X$, and therefore for all $x\in V$. But this 
contradicts Lemma 2.4.

\proclaim
Claim 4. The function $c_\flat:X\to\bbk$ is constant.
\endproclaim

We have to prove that $c_\flat(x_1)=c_\flat(x_2)$ whenever $x_1,x_2\in X$. If $x_2$ is 
a scalar multiple of $x_1$, then this equality is obviously true. Suppose that 
$x_1$ and $x_2$  are linearly independent. Note that
$$
c_\flat(x_1+x_2)\,\ell_{(x_1+x_2)y}=\ell'_{(x_1+x_2)y}
=\ell'_{x_1y}+\ell'_{x_2y}=c_\flat(x_1)\,\ell_{x_1y}+c_\flat(x_2)\,\ell_{x_2y}
$$
for all $y\in V$. If $\ell_{x_1y}$ and $\ell_{x_2y}$ are linearly independent 
for some $y$, then, comparing the coefficients, we deduce that both $c_\flat(x_1)$ 
and $c_\flat(x_2)$ are equal to $c_\flat(x_1+x_2)$, and we are done. We only have to 
be sure that such an element $y$ exists.

Otherwise we would have $\ell_{x_1y}\wedge\ell_{x_2y}=0$ for all $y\in V$. In 
this case $\ell_{x_2y}$ is a scalar multiple of $\ell_{x_1y}$ whenever 
$y\notin K_{x_1}$, and we can write $\ell_{x_2y}=\la(y)\ell_{x_1y}$ for some 
$\la(y)\in\bbk$. If $y_1,y_2\in V$ are linearly independent modulo the 
one-dimensional subspace $K_{x_1}$, then $\ell_{x_1y_1}$ and $\ell_{x_1y_2}$ 
are linearly independent, and, expressing $\ell'_{x_2(y_1+y_2)}$ as their 
linear combination, we deduce that $\la(y_1)=\la(y_1+y_2)=\la(y_2)$. It 
follows that $\la$ is a constant function.

In other words, there is $\la\in\bbk$ such that $\ell_{x_2y}=\la\ell_{x_1y}$ 
for all $y\in V\setm K_{x_1}$. By linearity in $y$ the last equality holds then 
for all $y\in V$, whence $\ell_{(x_2-\la x_1)y}=0$ for all $y\in V$. However, 
this contradicts Lemma 2.4.

\medskip
Thus Claim 4 is proved. It means that there is $c\in\bbk$ such that 
$\,\ell'_{xy}=c\,\ell_{xy}\,$ for all $x\in X$ and all $y\in V$. This equality 
holds then also for all $x\in\bbk X$, hence for all $x\in V$ by Claim 3. This 
shows that $\hbox{(v)}\Rar\hbox{(iv)}$.
\endproof

\section
3. The proof

We will verify condition (v) of Lemma 2.6. This will be done by analyzing 
identities (2.11) and (2.14) involving the linear functions $\ell_{xy}$ and 
$\ell'_{xy}$. Now one may safely forget the actual connection of these 
functions with Hecke symmetries. We will need only their properties stated in 
Lemma 2.4 and several simple consequences of (2.10).

Since the construction of $R$-symmetric algebras commutes with extensions of 
the base field, we may assume without loss of generality the field $\bbk$ to 
be algebraically closed. This will allow us to work comfortably with algebraic 
varieties and use some geometric arguments.

Denote by $U$ and $U'$ the subspaces of $V^*$ spanned by the sets 
$\{\ell_{xx}\mid x\in V\}$ and $\{\ell'_{xx}\mid x\in V\}$ respectively. 
Note that $U'=\{f\circ\ze\mid f\in U\}$ by (2.10). In particular, 
$\dim U'=\dim U$.

\proclaim
Proposition 3.1.
If $\,\dim U=1$, then $R$ commutes with $\ze\ot\ze$.
\endproclaim

\Proof.
Since $\dim U'=1$ too, we have $\dim(U+U')\le2$. Therefore $\bigwedge^2(U+U')$ 
is at most 1-dimensional subspace of $\bigwedge^2V^*$ which contains 
$\ell_{xx}\wedge\ell'_{yy}$ for any $x,y\in V$. Identity (2.11) yields
$$
\ell_{xy}\wedge\ell'_{xy}\in{\textstyle\bigwedge^2(U+U')}
$$
for all $x,y\in V$. The set $O=\{(x,y)\in V\times V\mid\ell_{xy}\notin U+U'\}$ 
is a Zariski open subset of $V\times V$. It is nonempty by Lemma 2.4. Note 
that for two linearly independent linear forms $\xi,\eta\in V^*$ the 
containment $\xi\wedge\eta\in\bigwedge^2(U+U')$ holds only when  both $\xi$ 
and $\eta$ lie in $U+U'$. Therefore for $(x,y)\in O$ the displayed 
containment implies that $\ell_{xy}\wedge\ell'_{xy}=0$.

On the other hand, 
$Z=\{(x,y)\in V\times V\mid\ell_{xy}\wedge\ell'_{xy}=0\}$ is a Zariski closed 
subset of an irreducible algebraic variety $V\times V$. We conclude that 
$Z=V\times V$ since $Z$ contains a nonempty Zariski open subset of $V\times V$. 
Thus $\ell_{xy}\wedge\ell'_{xy}=0$ for all $x,y\in V$. Now Lemma 2.6 applies.
\endproof

Proposition 3.1 confirms the conclusion of Theorem 1.1 in the case when $U$ 
has dimension 1. In fact this condition holds for all Hecke symmetries 
obtained by twisting from Hecke symmetries associated with the ordinary 
symmetric algebra $\bbS(V)$. It remains to consider the case $\,\dim U>1$ 
which could possibly lead to another class of Hecke symmetries.

\proclaim
Proposition 3.2.
Assume that $\,\chr\bbk\ne2$. There does not exist any Hecke symmetry $R$ 
such that the $R$-symmetric algebra\/ $\bbS(V,R)$ is isomorphic to some 
twisted polynomial algebra in $3$ indeterminates and for which $\,\dim U>1$.  
\endproclaim

The proof of Proposition 3.2 is much longer. It will be split into a series of 
lemmas which occupy the rest of this section. Let us assume that $\,\dim U>1$. 
We will see eventually that this assumption leads to a contradiction.

The next lemma provides a tool to derive certain relations of linear 
dependence between the linear forms $\ell_{xy}$ and $\ell'_{xy}$.

\proclaim
Lemma 3.3.
Let $X$ be an irreducible affine algebraic variety over an algebraically closed 
field $\bbk$ with a factorial coordinate ring $\bbk[X]$. 
Suppose that $W$ is a finite dimensional vector space over $\bbk$ and 
$\ph_1,\ldots,\ph_n,\psi:X\to W$ are morphisms in the category of algebraic 
varieties such that the set
$$
O=\{x\in X\mid\hbox{$\ph_1(x),\ldots,\ph_n(x)$ are linearly independent}\}
$$
is nonempty, while for each $x\in X$ the $n+1$ vectors 
$\ph_1(x),\ldots,\ph_n(x),\psi(x)\in W$ are linearly dependent. Then 
there exist functions $p_1,\ldots,p_n,p_0\in\bbk[X]$ such that 
$$
\gcd(p_1,\ldots,p_n,p_0)=1
\vadjust{\vskip-8pt}\eqno(3.1)
$$
and\vadjust{\vskip-8pt}
$$
p_0(x)\,\psi(x)=\sum_{i=1}^n\,p_i(x)\,\ph_i(x)
\quad\hbox{for all $x\in X$}.\eqno(3.2)
$$
Moreover, the equality $\,h_0(x)\,\psi(x)=\sum_{i=1}^nh_i(x)\,\ph_i(x)\,$ holds 
identically on $X$ for some collection of functions 
$\,h_1,\ldots,h_n,h_0\in\bbk[X]\,$ if and only if there is $\,g\in\bbk[X]\,$ 
such that $\,h_i=gp_i\,$ for all $i$.
\endproclaim

\Proof.
For each $x\in O$ the vector $\psi(x)$ is a linear combination of 
$\ph_1(x),\ldots,\ph_n(x)$ with uniquely determined coefficients. So
$$
\psi(x)=\sum_{i=1}^n\,f_i(x)\,\ph_i(x),\qquad x\in O,
$$
for some functions $f_i:O\to\bbk$. It is easy to see that $O$ is a Zariski 
open subset of $X$ and $f_1,\ldots,f_n$ are regular on $O$, hence rational on 
$X$. We can find $0\ne p_0\in\bbk[X]$ such that $p_0f_i=p_i\in\bbk[X]$ for 
each $i=1,\ldots,n$. Equality (3.2) is satisfied on some nonempty Zariski open 
subset of $X$, but then it must hold everywhere. Factoring out the 
greatest common divisor of $p_1,\ldots,p_n,p_0$ we achieve (3.1).

In the second assertion we have
$$
\sum_{i=1}^n\,\bigl(p_0(x)\,h_i(x)-h_0(x)\,p_i(x)\bigr)\,\ph_i(x)=0\quad
\hbox{for all $x\in X$}.
$$
It follows from the definition of the set $O$ that 
$p_0(x)h_i(x)=h_0(x)p_i(x)$ for all $x\in O$ and $i=1,\ldots,n$. Any two 
regular functions $X\to\bbk$ must be equal provided that they agree on a 
nonempty Zariski open subset of $X$. Hence $p_0h_i=h_0p_i$ for each $i$. 
Factoriality of $\bbk[X]$ together with (3.1) ensures that $p_0$ divides $h_0$ 
in the ring $\bbk[X]$. So $\,g=h_0p_0^{-1}$ will do.
\endproof

We will apply Lemma 3.3 in the situation where $X$ is either $V$ or 
$V\times V$. So $X$ is a vector space, and $\bbk[X]\cong\bbS(X^*)$ is the 
algebra of polynomial functions on $X$ generated by the dual space $X^*$ of 
linear functions on $X$.

Polynomial functions on $V\times V$ are functions of two arguments taken in $V$. 
We will encounter polynomial functions $f:V\times V\to\bbk$ which are 
homogeneous in each of the two arguments. We say that $(m,n)$ is the 
\emph{bidegree} of $f$ if $f$ is homogeneous of degree $m$ 
in the first argument and homogeneous of degree $n$ in the second. The 
\emph{total degree} of $f$ is the sum $m+n$ of its bidegree components. 
Note that any divisor of a bihomogeneous polynomial function is itself 
bihomogeneous.

Each of the expressions $\ell_{xy}$, $\ell'_{xy}$, $\ell_{xx}$, $\ell'_{yy}$ 
regarded as a function of the pair $(x,y)$ gives a quadratic polynomial map 
$V\times V\to V^*$. Since $\dim U>1$, there exist $x,y\in V$ such that 
$\ell_{xx}$ and $\ell'_{yy}$ are linearly independent. On the other hand, the 
three elements $\ell_{xy}$, $\ell_{xx}$, $\ell'_{yy}$ of the vector space 
$V^*$ are always linearly dependent, and so too are 
$\ell'_{xy}$, $\ell_{xx}$, $\ell'_{yy}$ since 
$$
\ell_{xy}\wedge\ell_{xx}\wedge\ell'_{yy}=0\qquad\hbox{and}\qquad
\ell'_{xy}\wedge\ell_{xx}\wedge\ell'_{yy}=0\eqno(3.3)
$$
for all $x,y\in V$ by (2.11). Lemma 3.3 provides polynomial 
functions $p_i,p'_i\in\bbk[V\times V]$ for $i=0,1,2$ such that
$$
\gcd(p_0,p_1,p_2)=\gcd(p'_0,p'_1,p'_2)=1,\eqno(3.4)
$$
and
$$
\openup1\jot
\eqalign{
p_0(x,y)\,\ell_{xy}&=p_1(x,y)\,\ell_{xx}+p_2(x,y)\,\ell'_{yy}\,,\cr
p'_0(x,y)\,\ell'_{xy}&=p'_1(x,y)\,\ell_{xx}+p'_2(x,y)\,\ell'_{yy}
}\eqno(3.5)
$$
for all $x,y\in V$. As seen from the final assertion of Lemma 3.3, those $p_i$ 
and $p'_i$ are polynomial functions of the smallest degree that can occur in 
such relations of linear dependence. Since every polynomial function on 
$V\times V$ is a sum of its bihomogeneous components, the functions $p_i$ and 
$p'_i$ are necessarily bihomogeneous.

We will need very precise information about the functions $p_i$, $p'_i$. As a 
first step we deduce an identity which connects them. Making use of (3.5), we get
$$
p_0(x,y)\,p'_0(x,y)\,\ell_{xy}\wedge\ell'_{xy}
=\bigl(p_1(x,y)\,p'_2(x,y)-p_2(x,y)\,p'_1(x,y)\bigr)\,\ell_{xx}\wedge\ell'_{yy}.
$$
The rule $(x,y)\mapsto\ell_{xx}\wedge\ell'_{yy}$ defines a nonzero polynomial 
map $V\times V\to\bigwedge^2V^*$. Since 
$\,\ell_{xy}\wedge\ell'_{xy}=\ell_{xx}\wedge\ell'_{yy}\,$ and the ring 
$\bbk[V\times V]$ is a domain, it follows that
$$
p_0p'_0=p_1p'_2-p_2p'_1.
$$
We will see later that $p'_0$ can be taken equal to $p_0$. In this case the 
previous identity is written as
$$
p_0^2=p_1p'_2-p_2p'_1.\eqno(3.6)
$$

Take any basis $\ell_1,\ell_2,\ell_3$ for the vector space 
$V^*$. We can write
$$
\eqalign{
\ell_{xx}&=\sum a_i(x)\,\ell_i,\qquad\ell_{xy}=\sum b_i(x,y)\,\ell_i,\cr
\ell'_{xx}&=\sum a'_i(x)\,\ell_i,\qquad\ell'_{xy}=\sum b'_i(x,y)\,\ell_i
}\eqno(3.7)
$$
for some quadratic forms $a_i,a'_i$ and bilinear forms $b_i,b'_i$ such that 
$a_i(x)=b_i(x,x)$ and $a'_i(x)=b'_i(x,x)$. Now
$$
\ell_{xx}\wedge\ell'_{yy}=
\De_1(x,y)\,\ell_2\wedge\ell_3+\De_2(x,y)\,\ell_1\wedge\ell_3+
\De_3(x,y)\,\ell_1\wedge\ell_2\eqno(3.8)
$$
where $\De_i(x,y)$, $i=1,2,3$, are the minors of order 2 of the matrix
$$
\pmatrix{a_1(x)&a_2(x)&a_3(x)\cr a'_1(y)&a'_2(y)&a'_3(y)}.
$$

\proclaim
Lemma 3.4.
For each $i=1,2,3$ the polynomial function $\De_i$ is divisible by $p_0$ and 
by $p'_0$ in the ring $\,\bbk[V\times V]$.
\endproclaim

\Proof.
In terms of coordinate representation (3.7) the first equation in (3.3) means 
that the matrix
$$
\pmatrix{b_1(x,y)&b_2(x,y)&b_3(x,y)\cr
a_1(x)&a_2(x)&a_3(x)\cr
a'_1(y)&a'_2(y)&a'_3(y)}
$$
has identically zero determinant. Hence there are 3 different relations of 
linear dependence between the rows of this matrix, the coefficients being the 
second order minors each time extracted from some pair of columns of the matrix. 
This gives 3 relations of linear dependence between $\ell_{xy}$, $\ell_{xx}$, 
$\ell'_{yy}$ in which the respective coefficient of $\ell_{xy}$ is 
$\De_i(x,y)$ for $i=1,2,3$. For example,
$$
\left|\matrix{a_1(x)&a_2(x)\cr a'_1(y)&a'_2(y)}\right|\ell_{xy}
-\left|\matrix{b_1(x,y)&b_2(x,y)\cr a'_1(y)&a'_2(y)}\right|\ell_{xx}
+\left|\matrix{b_1(x,y)&b_2(x,y)\cr a_1(x)&a_2(x)}\right|\ell'_{yy}=0.
$$
The final assertion in Lemma 3.3 shows that $p_0$ divides $\De_i$. Working with 
the second equation in (3.3) we deduce similarly that $p'_0$ divides $\De_i$.
\endproof

Consider the greatest common divisor
$$
d=\gcd(\De_1,\De_2,\De_3)\in\bbk[V\times V].\eqno(3.9)
$$
Since $\De_1,\De_2,\De_3$ are bihomogeneous of bidegree $(2,2)$, the 
function $d$ has to be bihomogeneous too, and its degree in each of the two 
arguments cannot exceed 2. The same can be said about $p_0$ and $p'_0$ which 
are divisors of $d$ by Lemma 3.4. 

\proclaim
Lemma 3.5.
The three polynomial functions $\De_1,\De_2,\De_3$ are linearly independent. 
As a consequence, $\,\deg d<4$.
\endproclaim

\Proof.
Suppose that $\De_1,\De_2,\De_3$ are linearly dependent. It is then seen from 
(3.8) that the set $\{\ell_{xx}\wedge\ell'_{yy}\mid x,y\in V\}$ spans a proper 
subspace of $\bigwedge^2V^*$. However, the linear span of this set is a 
subspace containing $\xi\wedge\eta$ for all $\xi\in U$ and $\eta\in U'$, and 
we claim that it is the whole space $\bigwedge^2V^*$.

If $U=V^*$, then $U'=V^*$ as well, and the claim is obvious. Suppose that 
$U\ne V^*$. Then $\dim U=\dim U'=2$. If $U\ne U'\!$, then the claim is still 
true. Finally, if $U'=U$, then the argument in the proof of Proposition 3.1 
shows that actually $\ell_{xy}\wedge\ell'_{xy}=0$ for all $x,y\in V$. By 
(2.11) in this case $\ell_{xx}\wedge\ell'_{yy}=0$ for all $x,y\in V$, but this 
contradicts the assumption $\,\dim U=2$.
\endproof

\proclaim
Lemma 3.6.
The polynomial functions $p_i,$ $p'_i$ are all nonzero. As a consequence, $p_0$ 
and $p'_0$ are not functions of just one argument. In particular, $\deg p_0>1$ 
and $\deg p'_0>1$.
\endproclaim

\Proof.
Suppose that $p_2=0$, for example. By (3.5) and (3.7) then
$$
p_0(x,y)\,b_i(x,y)=p_1(x,y)\,a_i(x)
$$
for all $x,y\in V$ and each $i=1,2,3$. Since the ring $\bbk[V\times V]$ is 
factorial and since $\gcd(p_0,p_1)=1$ by (3.4), the function $p_1$ must divide 
each $b_i$ in that ring.

The functions $b_1,b_2,b_3$ have bidegree $(1,1)$. Hence each bidegree component 
of $p_1$ does not exceed 1. If $p_1$ has bidegree $(1,1)$, then each $b_i$ will 
be a scalar multiple of $p_1$, and so $\ell_{xy}=p_1(x,y)\,\ell_0$ for some 
linear form $\ell_0\in V^*$ which depends neither on $x$ nor on $y$. In this 
case all linear forms $\ell_{xy}$ are contained in the one-dimensional 
subspace of $V^*$ spanned by $\ell_0$, in contradiction with Lemma 2.4. 
Clearly $p_1$ must depend on $y$. Therefore the only possibility left is that 
$p_1$ has bidegree $(0,1)$, i.e., $p_1(x,y)=\xi(y)$ for some $\xi\in V^*$. But 
then $\ell_{xy}=0$ for all $x\in V$ and $y\in\Ker\xi$, again in contradiction 
with Lemma 2.4.

It is proved similarly that $p_1\ne0$, while $p_0\ne0$ is clear already from 
the assumption that $\ell_{xx}\wedge\ell'_{yy}$ does not vanish identically. 

Let $(m_i,n_i)$ be the bidegree of $p_i$. By the uniqueness of expressions 
(3.5) the basis linear forms $\ell_1,\ell_2,\ell_3$ must occur with 
coefficients of equal bidegrees in each term of the first equation there. 
It follows that
$$
m_0=m_1+1=m_2-1\qquad\hbox{and}\qquad n_0=n_1-1=n_2+1.\eqno(3.10)
$$
Hence $m_0>0$ and $n_0>0$. The same argument applies to $p'_0$, $p'_1$, $p'_2$.
\endproof

\proclaim
Lemma 3.7.
The equality $\,\ell_{xx}\wedge\ell'_{xx}=0\,$ holds for all $x\in V$. As a 
consequence, $U=U'=V^*$.
\endproclaim

\Proof.
Setting $y=x$ in (3.5), we get
$$
\eqalign{
\bigl(p_0(x,x)-p_1(x,x)\bigr)\,\ell_{xx}&=p_2(x,x)\,\ell'_{xx}\,,\cr
\bigl(p'_0(x,x)-p'_2(x,x)\bigr)\,\ell'_{xx}&=p'_1(x,x)\,\ell_{xx}\,.
}\eqno(3.11)
$$
If $p_0$ has degree 1 in the second argument, then by (3.10) $p_2$ must have 
degree 0 in the second argument, i.e., $p_2$ depends only on the first 
argument. In this case $p_2$ is not identically zero on the diagonal 
$$
D=\{(x,x)\mid x\in V\}\sbs V\times V,
$$
and for all $x\in V$ such that $p_2(x,x)\ne0$ the linear function $\ell'_{xx}$ 
is a scalar multiple of $\ell_{xx}$ in view of (3.11). Similarly, if $p'_0$ 
has degree 1 in the first argument, then $p'_1$ depends only on the second 
argument. In this case $p'_1$ is not identically zero on $D$, and so 
$\ell_{xx}$ is a scalar multiple of $\ell'_{xx}$ for all $x\in V$ such that 
$p'_1(x,x)\ne0$. In both cases the equality $\,\ell_{xx}\wedge\ell'_{xx}=0\,$ 
holds for all $x$ in a nonempty Zariski open subset of $V$, and therefore 
everywhere on $V$.

To confirm the first assertion of Lemma 3.7 it remains to note that one of the 
previous two conditions is satisfied in any event. Since $p_0$ and $p'_0$ are 
divisors of $d$, it follows from Lemma 3.6 that $d$ has positive degree in 
each of the two arguments. But $\deg d<4$ by Lemma 3.5. This leaves only 3 
possibilities for the bidegree of $d$. It is $(1,1)$, or $(1,2)$, or $(2,1)$. 
If neither $p_0$ nor $p'_0$ have bidegree $(1,1)$, then both functions will be 
scalar multiples of $d$, and so they will have the same bidegree, either 
$(1,2)$ or $(2,1)$.

The equality $\,\ell_{xx}\wedge\ell'_{xx}=0\,$ implies that $\ell'_{xx}$ is a 
scalar multiple of $\ell_{xx}$, and so $\ell'_{xx}\in U$, whenever 
$\ell_{xx}\ne0$. Thus $\ell'_{xx}\in U$ for all $x$ in a nonempty Zariski open 
subset of $V$, and therefore for all $x\in V$. It follows that $U'\sbs U$, and 
in fact $U'=U$ since these two spaces have the same dimension. Now 
$$
p_0(x,y)\,\ell_{xy}\in U+U'=U
$$
for all $x,y\in V$ according to (3.5). Hence $\ell_{xy}\in U$ whenever 
$p_0(x,y)\ne0$. In other words, $\ell_{xy}\in U$ for all pairs $(x,y)$ in a 
nonempty Zariski open subset of $V\times V$, but then this containment holds 
everywhere on $V\times V$, and the equality $U=V^*$ follows from Lemma 2.4.
\endproof

\proclaim
Lemma 3.8.
The assumption that the three quadratic forms $a_1,a_2,a_3$ have a nonscalar 
common divisor in the ring $\bbk[V]$ leads to a contradiction. Hence we must 
have $\,\gcd(a_1,a_2,a_3)=1$. Also, $\,\gcd(a'_1,a'_2,a'_3)=1$.
\endproclaim

\Proof.
Suppose that $\xi\in\bbk[V]$ is a nonscalar polynomial function which divides 
each $a_i$. Clearly, $\xi$ is homogeneous of degree $\le2$. If $\deg\xi=2$ 
then each $a_i$ is a scalar multiple of $\xi$, and so 
$\ell_{xx}=\xi(x)\ell_0$ for some linear function $\ell_0\in V^*$ which does 
not depend on $x$. This contradicts the assumption $\dim U>1$, however.

Therefore $\deg\xi=1$, i.e., $\xi\in V^*$. We have $\ell_{xx}=0$ for all 
$x\in\Ker\xi$, whence $\ell'_{xx}=0$ for all $x\in\ze^{-1}(\Ker\xi)$ in view of 
(2.10). It follows that each $a'_i$ is divisible in $\bbk[V]$ by the linear 
function $\eta=\xi\circ\ze$. We can write 
$$ 
\ell_{xx}=\xi(x)\sum f_i(x)\ell_i\,,
\qquad\ell'_{xx}=\eta(x)\sum g_i(x)\ell_i\eqno(3.12)
$$
for some linear forms $f_i,g_i\in V^*$, and we get 
$\,\De_i(x,y)=\xi(x)\,\eta(y)\,\de_i(x,y)\,$ where $\de_i(x,y)$ is the 
respective minor of order 2 of the matrix
$$
\pmatrix{f_1(x)&f_2(x)&f_3(x)\cr g_1(y)&g_2(y)&g_3(y)}.
$$
Each $\de_i$ is a polynomial function on $V\times V$ of bidegree $(1,1)$, and 
$\de_1,\de_2,\de_3$ are linearly independent by Lemma 3.5. Hence the degree of 
any common divisor of these three functions is less than 2.

Suppose that $\gcd(\de_1,\de_2,\de_3)\ne1$. Then there exists a linear form 
$h$ on $V\times V$ of bidegree either $(1,0)$ or $(0,1)$ which divides each 
$\de_i$. If, for example, the bidegree of $h$ is $(1,0)$, then 
$\Ker h=L\times V$ where $L$ is a subspace of codimension 1 in $V$, and each 
$\de_i$ must vanish on $L\times V$. This means that
$$
f_j(x)\,g_k(y)=f_k(x)\,g_j(y)\quad\hbox{for all $\,j,k\in\{1,2,3\}$, 
$\,x\in L$, $\,y\in V$}.
$$
Since $\dim U'=3$ by Lemma 3.7, the functions $g_1,g_2,g_3$ must be linearly 
independent, and it follows from the displayed equalities above that 
each $f_i$ vanishes on $L$. But then $f_1,f_2,f_3$ are scalar multiples of 
one function, and we see from (3.12) that all functions $\ell_{xx}$, 
$\,x\in V$, lie in a one-dimensional subspace of $V^*$. This contradicts the 
assumption $\dim U>1$, however.

We conclude that $\gcd(\de_1,\de_2,\de_3)=1$, and therefore 
$\gcd(\De_1,\De_2,\De_3)=d$ where $d$ is defined by the formula 
$$
d(x,y)=\xi(x)\,\eta(y),\qquad x,y\in V.
$$
Since $p_0$ and $p'_0$ are divisors of $d$ of degree at least $2$ by Lemma 
3.6, both $p_0$ and $p'_0$ must be scalar multiples of $d$.
By scaling the functions we may assume without loss of generality that 
$$
p_0=p'_0=d.
$$
Here $p_0$ has bidegree $(1,1)$. This implies that $p_1$ and $p'_1$ have 
bidegree $(0,2)$, while $p_2$ and $p'_2$ have bidegree $(2,0)$.

Multiplying the first equation in (3.5) by $p'_2(x,y)$, subtracting then the 
second equation multiplied by $p_2(x,y)$, and making use of (3.6), we get
$$
\openup1\jot
\eqalign{
p_0(x,y)\,\bigl(p'_2(x,y)\,\ell_{xy}-p_2(x,y)\,\ell'_{xy}\bigr)
&=\bigl(p'_2(x,y)\,p_1(x,y)-p_2(x,y)\,p'_1(x,y)\bigr)\,\ell_{xx}\cr
&=p_0(x,y)^2\,\ell_{xx}
}
$$
for all $x,y\in V$. Since the ring $\bbk[V\times V]$ is a domain, the common 
factor $p_0(x,y)$ in the left and right hand sides can be cancelled out, which 
results in the identity
$$
p'_2(x,y)\,\ell_{xy}-p_2(x,y)\,\ell'_{xy}=p_0(x,y)\,\ell_{xx}\,.
$$
Hence
$$
p'_2(x,y)\,\ell_{xy}(x)-p_2(x,y)\,\ell'_{xy}(x)=p_0(x,y)\,\ell_{xx}(x)
=\xi(x)\,\eta(y)\,\ell_{xx}(x)\,.
$$
We have $p'_2(x,y)=p'_2(x,0)$ and $p_2(x,y)=p_2(x,0)$ since $p_2$ and $p'_2$ 
depend only on the first argument. Now making use of (2.14), we rewrite the 
preceding equality as
$$
p'_2(x,0)\,\ell_{xx}(y)-p_2(x,0)\,\ell'_{xx}(y)=\xi(x)\,\ell_{xx}(x)\,\eta(y)\,.
$$
Since this holds for all $x,y\in V$, we get
$$
p'_2(x,0)\,\ell_{xx}-p_2(x,0)\,\ell'_{xx}=\xi(x)\,\ell_{xx}(x)\,\eta\eqno(3.13)
$$
for all $x\in V$. The subset
$$
O=\{x\in V\mid\hbox{$\ell_{xx}$ is not a scalar multiple of $\eta$}\}
$$
is clearly Zariski open in $V$. It is also nonempty in view of Lemma 3.7. 
Suppose that $x\in O$. Then $\ell_{xx}\ne0$, and the linear form in the left 
hand side of (3.13) is a scalar multiple of $\ell_{xx}$ since so is 
$\ell'_{xx}$ by Lemma 3.7. Since $\ell_{xx}$ and $\eta$ are linearly 
independent, both sides of (3.13) must vanish. This yields 
$\,\xi(x)\,\ell_{xx}(x)=0$.

We conclude that $\,\xi(x)\,\ell_{xx}(x)=0\,$ for all $x\in V$ since this equality 
holds for all $x$ in a nonempty Zariski open subset of $V$. Since $\xi\ne0$, 
we get
$$
\ell_{xx}(x)=0\quad\hbox{for all $x\in V$}.
$$
With this at hand we go back to (3.12). Consider the bilinear form 
$\ga:V\times V\to\bbk$ and linear forms $\ga_x\in V^*$ defined by the rule
$$
\ga(x,y)=\sum f_i(x)\,\ell_i(y),\qquad\ga_x(y)=\ga(x,y)
$$
for $x,y\in V$. Then $\,\ell_{xx}=\xi(x)\,\ga_x$. Hence 
$\,\xi(x)\,\ga(x,x)=\ell_{xx}(x)=0\,$ for all $x\in V$, and therefore $\ga(x,x)=0$ 
for all $x\in V$. This shows that the bilinear form $\ga$ is alternating, 
which implies that its rank is even. Since $\dim V=3$, any alternating 
bilinear form on $V$ is degenerate. Hence there exists $0\ne v\in V$ such that 
$\ga_x(v)=0$, and therefore also $\ell_{xx}(v)=0$, for all $x\in V$. This 
contradicts the equality $U=V^*$ of Lemma 3.7.

At the beginning of the proof we have seen that the assumption 
$\,\gcd(a_1,a_2,a_3)\ne1$ implies that $\,\gcd(a'_1,a'_2,a'_3)\ne1$. The 
opposite implication is proved quite similarly. Therefore the assumption 
$\,\gcd(a'_1,a'_2,a'_3)\ne1$ leads to a contradiction too.
\endproof

\proclaim
Lemma 3.9.
The polynomial function $d=\gcd(\De_1,\De_2,\De_3)$ has bidegree $(1,1)$, i.e., 
$d:V\times V\to\bbk$ is a bilinear form on $V$. As a consequence, both $p_0$ 
and $p'_0$ are scalar multiples of $d$.
\endproclaim

\Proof.
For each $x\in V$ the linear forms $\ell_{xx}$ and $\ell'_{xx}$ are linearly 
dependent by Lemma 3.7. Applying Lemma 3.3 to the two polynomial maps 
$V\to V^*$ given by the assignments $x\mapsto\ell_{xx}$ and $x\mapsto\ell'_{xx}$, 
respectively, we deduce that there are polynomial functions 
$\,\xi,\eta\in\bbk[V]\,$ such that $\,\gcd(\xi,\eta)=1\,$ and
$$
\xi(x)\,\ell'_{xx}=\eta(x)\,\ell_{xx}\quad\hbox{for all $x\in V$}.
$$
In the ring $\bbk[V]$ we have then $\xi a'_i=\eta a_i$ for each $i=1,2,3$ (see 
(3.7)). Hence $\xi$ is a common divisor of quadratic forms $a_1,a_2,a_3$, 
and so $\xi\in\bbk$ by Lemma 3.8. Clearly, $\deg\xi=\deg\eta$, whence 
$\eta\in\bbk$ as well. It follows that there is $0\ne c\in\bbk$ such that 
$a'_i=ca_i$ for each $i$. But then each $\De_i$ is skewsymmetric, i.e., 
$$
\De_i(y,x)=-\De_i(x,y)\quad\hbox{for all $x,y\in V$},
$$
and it follows that the degree of $d$ in the second argument must be the same 
as its degree in the first argument. So the total degree of $d$ is even and 
less than 4 by Lemma 3.5. Since $p_0$ and $p'_0$ are divisors of $d$, there is 
no possibility other than
$$
\deg d=\deg p_0=\deg p'_0=2,
$$
and the bidegree of these functions is necessarily $(1,1)$.
\endproof

Thus we may take $p_0=p'_0=d$ in (3.5).

\proclaim
Lemma 3.10.
The function $p_0$ is a bilinear form of rank $1$, i.e., $p_0$ is the product 
of two linear forms of bidegrees $(1,0)$ and $(0,1)$.
\endproclaim

\Proof.
Identity (3.6) shows that $p_0$ is contained in the radical of the ideal of 
the ring $\bbk[V\times V]$ generated by $p_1$ and $p_2$. Hence $p_0$ vanishes 
on the set of common zeros of the two functions $p_1$ and $p_2$. Note that 
$p_1$ and $p'_1$ have bidegree $(0,2)$, while $p_2$ and $p'_2$ have bidegree 
$(2,0)$. Hence the zero sets of $p_1$ and $p_2$ are, respectively, $V\times Y$ 
and $X\times V$ where $X$ and $Y$ are some quadratic conical hypersurfaces 
in $V$. The common zero set of $p_1$ and $p_2$ is $X\times Y$. Since $p_0$ is 
bilinear, it vanishes on $\bbk X\times\bbk Y$ where $\bbk X$ and $\bbk Y$ are 
the linear spans of $X$ and $Y$.

Since $p_0\ne0$, the equalities $\bbk X=V$ and $\bbk Y=V$ cannot hold 
simultaneously. Therefore either $X$ or $Y$ must be a linear subspace of 
codimension 1 in $V$. This means that at least one of the two functions 
$p_1$, $p_2$ is the square of a linear form on $V\times V$. If $\bbk X=V$ then 
$Y$ is a linear subspace of $V$ contained in the right kernel of $p_0$, whence 
$\rk p_0=1$. Similarly, $\rk p_0=1$ whenever $\bbk Y=V$.

Suppose that $\rk p_0\ne1$. Then both $X$ and $Y$ are linear subspaces of 
codimension 1 in $V$. Hence there exist $\xi,\eta\in V^*$ such that
$$
p_1(x,y)=\eta(y)^2,\qquad p_2(x,y)=\xi(x)^2
$$
for all $x,y\in V$. A similar argument shows that there are 
$\xi',\eta'\in V^*$ such that
$$
p'_1(x,y)=\eta'(y)^2,\qquad p'_2(x,y)=\xi'(x)^2
$$
for all $x,y\in V$. Identity (3.6) is now written as
$$
p_0(x,y)^2=\eta(y)^2\,\xi'(x)^2-\xi(x)^2\,\eta'(y)^2.
$$
Hence $p_0^2=fg$ where $f,g\in\bbk[V\times V]$ are defined by the formulas
$$
f(x,y)=\eta(y)\,\xi'(x)-\xi(x)\,\eta'(y),\qquad
g(x,y)=\eta(y)\,\xi'(x)+\xi(x)\,\eta'(y).
$$
Since $p_0$ is a bilinear form of rank $>1$, it is an irreducible element of 
the factorial ring $\bbk[V\times V]$, and it follows from the factoriality 
that $f$ and $g$ are both scalar multiples of $p_0$. Then so too is $f+g$. But
$$
(f+g)(x,y)=2\,\eta(y)\,\xi'(x),
$$
in contradiction with irreducibility of $p_0$. Note that here we do need the 
assumption that $\,\chr\bbk\ne2$.
\endproof

\proclaim
Lemma 3.11.
There are $\xi,\eta\in V^*$ such that
$$
p_0(x,y)=\xi(x)\,\eta(y)\eqno(3.14)
$$
and either $\,p_1(x,y)=\eta(y)^2\,$ or $\,p_2(x,y)=\xi(x)^2\,$ for all $x,y\in V$.
\endproclaim

\Proof.
We continue the argument in the proof of Lemma 3.10. Since $p_0$ has rank 1, 
its left kernel $L=\{x\in V\mid p_0(x,V)=0\}$ is a vector subspace of 
codimension 1 in $V$. If $X=L$, then, as we have seen in Lemma 3.10, $p_2$ is 
a square, i.e., $p_2=\xi_1^{\,2}$ where $\xi_1$ is a linear form of 
bidegree $(1,0)$. There is $\xi\in V^*$ such that $\,\xi_1(x,y)=\xi(x)\,$ 
for all $x,y\in V$. Since the function $p_0$ vanishes on 
$X\times V=\Ker\xi_1$, it is divisible by $\xi_1$ in the ring 
$\bbk[V\times V]$. This yields factorization (3.14).

Suppose that $X\ne L$. Then $X\not\sbs L$ since $X$ is a subvariety of 
codimension 1 in $V$, and it follows that $\bbk X+L=V$. Since $p_0$ vanishes 
on $(\bbk X+L)\times Y$, its right kernel $K=\{y\in V\mid p_0(V,y)=0\}$ 
contains $Y$. Moreover, $K=Y$ since $Y$ is a subvariety of codimension 1 in 
$V$. In particular, $Y$ is a linear subspace. This implies that 
$p_1=\eta_2^{\,2}$ where $\eta_2$ is a linear form of bidegree 
$(0,1)$. There is $\eta\in V^*$ such that $\,\eta_2(x,y)=\eta(y)\,$ for 
all $x,y\in V$. Since $p_0$ vanishes on $V\times Y=\Ker\eta_2$, it is 
divisible by $\eta_2$, whence formula (3.14) with a suitable $\xi$.
\endproof

The rest of the proof will go really fast. In the first equation of (3.5) we 
compare the coefficients of each basis linear form $\ell_i$. This gives
$$
p_0(x,y)\,b_i(x,y)=p_1(x,y)\,a_i(x)+p_2(x,y)\,a'_i(y)
$$
for all $x,y\in V$. There are two possibilities described in Lemma 3.11. In 
one case both $p_0$ and $p_2$ are divisible by the linear form $\xi_1$ of 
bidegree $(1,0)$ corresponding to $\xi$. On the other hand, $\xi_1$ cannot 
divide the function $p_1$ of bidegree $(0,2)$, and it follows from the 
displayed equality above and factoriality that each $a_i$ is divisible by 
$\xi$. In another case both $p_0$ and $p_1$ are divisible by the linear form 
$\eta_2$ of bidegree $(0,1)$ corresponding to $\eta$, but $\eta_2$ 
does not divide $p_2$. This implies that each $a'_i$ is divisible by $\eta$. 
Now Lemma 3.8 eliminates both possibilities.
\endproof

This completes the proof of Proposition 3.2 and Theorem 1.1.

\section
4. Explicit determination of Hecke symmetries

As in Theorem 1.1 we assume in this section that $\chr\bbk\ne2$. Moreover, 
everywhere with the exception of Proposition 4.1 and preceding definitions it 
will be assumed that $\bbk$ is algebraically closed (actually it suffices to 
assume that $\bbk$ contains square roots of all its elements). This is 
required by the use of explicit formulas obtained in \cite{Sk23b}.

Given a graded factor algebra $A$ of the tensor algebra $\bbT(V)$, we denote by
$$
\HeckeSym(A)
$$
the set of all Hecke symmetries on the vector space $V$ such 
that $\bbS(V,R)=A$ where the exact equality means that the two algebras are 
factor algebras of $\bbT(V)$ by the same ideal. In the case when $A=\bbT(V)/I$ 
where $I$ is the ideal of $\bbT(V)$ generated by 
$\{\ze(x)y-\ze(y)x\mid x,y\in V\}\sbs V^{\ot2}$ for $\ze\in\GL(V)$ we write 
$\bbS(V)_\ze$ instead of $A$, thus identifying the twisted algebra 
$\bbS(V)_\ze$ with a factor algebra of $\bbT(V)$ by means of an isomorphism of 
graded algebras which acts as the identity operator on homogeneous elements of 
degree 1.

By Theorem 1.1 the set $\HeckeSym(\,\bbS(V)_\ze)$ consists of the $\ze$-twists 
$R_\ze$ of those Hecke symmetries $R$ in the set $\HeckeSym(\,\bbS(V))$ which 
commute with $\ze\ot\ze$. By \cite{Sk23b, Theorem 5.1} each $R$ in the latter 
set is given by the formula
$$
R(xy)={q-1\over2}xy+{q+1\over2}yx-g(x,y)\,a\wdg b-x\wdg Ty-y\wdg Tx,\quad 
x,y\in V,\eqno(4.1)
$$
where $a,b\in V$ are two vectors, $g:V\times V\to\bbk$ a symmetric bilinear 
form satisfying
$$
(q-1)^2=4\,\bigl(g(a,b)^2-g(a,a)\,g(b,b)\bigr)\eqno(4.2)
$$
and $T:V\to V$ the linear operator defined by the rule
$$
Tx=g(b,x)\,a-g(a,x)\,b,\qquad x\in V.\eqno(4.3)
$$
The product $\wdg$ used in (4.1) is the one defined in section 1, however with 
respect to the identity operator, i.e., $x\wdg y=xy-yx$ for $x,y\in V$.

The linear operator given by (4.1) depends on the pair $(t,g)$ where 
$t=a\wedge b$ is a bivector. A different choice of vectors 
$a,b$ making the same bivector does not change the operator. The parameter $q$ 
is determined by relation (4.2) in which the right hand side is a function of 
the pair $(t,g)$. If either $t=0$ or $g=0$, then this operator is the flip 
$R_0$ which sends $xy$ to $yx$ for all $x,y\in V$. If $t\ne0$ and $g\ne0$, 
then another pair $(t',g')$ produces the same operator if and only if $t'=ct$ 
and $g'=c^{-1}g$ for some $0\ne c\in\bbk$.

This leads to the following parametrization of the set $\HeckeSym(\,\bbS(V))$. 
Define
$$
\De(a\wedge b,\,g)=g(a,a)\,g(b,b)-g(a,b)^2\eqno(4.4)
$$
and denote by $P$ the set of all triples $(t,g,q)$ where $t\in\bigwedge^2V$ is 
a nonzero bivector, $g:V\times V\to\bbk$ a nonzero symmetric bilinear form, 
$q\in\bbk$ a nonzero scalar such that
$$
(q-1)^2=-4\,\De(t,g).\eqno(4.5)
$$
The multiplicative group $\bbk^\times$ of the field $\bbk$ acts on $P$ 
according to the rule
$$
c\cdot(t,\,g,\,q)=(c{\mskip1mu}t,\,c^{-1}g,\,q),\qquad c\in\bbk^\times.\eqno(4.6)  
$$
The elements of the set $\HeckeSym(\,\bbS(V))\setm\{R_0\}$ are then in a 
bijective correspondence with the $\bbk^\times$-orbits in $P$. Furthermore, 
$\bigwedge^2V$ and the space of symmetric bilinear forms on $V$ are 
$\GL(V)$-modules in a natural way. These module structures give rise to an 
action of $\GL(V)$ on $P$ under which $q$ remains unaffected. Conjugation by 
the operators $\ph\ot\ph$, $\,\ph\in\GL(V)$, is the corresponding action of 
$\GL(V)$ on Hecke symmetries. Therefore two Hecke symmetries in the set 
$\,\HeckeSym(\,\bbS(V))\setm\{R_0\}\,$ are equivalent if and only if they 
correspond to two triples in the set $P$ lying in the same 
$(\GL(V)\times\bbk^\times)$-orbit.

The flip $R_0$ commutes with $\ze\ot\ze$ for any $\ze\in\GL(V)$. The 
corresponding twist ${R_0}_\ze$ is the linear operator on $V^{\ot2}$ such that
$$
{R_0}_\ze(xy)=\ze(y)\,\ze^{-1}(x)\quad\hbox{for $x,y\in V$}.\eqno(4.7)
$$
It may be viewed as a distinguished element of the set 
$\,\HeckeSym(\,\bbS(V)_\ze)$. The remaining Hecke symmetries in this set can 
be described in terms of triples in $P$:

\proclaim
Proposition 4.1.
For $\ze\in\GL(V)$ put
$$
\openup1\jot
\displaylines{
P(\ze)=\{(t,g,q)\in P\mid\,\ze\cdot(t,g,q)\in\,\bbk^\times\!\cdot(t,g,q)\},\cr
G(\ze)=\{\,\ph\in\GL(V)\mid
\,\hbox{$\ph\mskip1mu\ze\ph^{-1}$ is a scalar multiple of $\ze$}\}.
}
$$
There is a bijection between the set $\,\HeckeSym(\,\bbS(V)_\ze)\setm\{{R_0}_\ze\}$ 
and the set $P(\ze)/\mskip1mu\bbk^\times$ of $\,\bbk^\times$-orbits in $P(\ze)$. 
Under this bijection the equivalence classes of Hecke symmetries correspond
to the $G(\ze)$-orbits in $P(\ze)/\mskip1mu\bbk^\times$.
\endproclaim

\Proof.
If $R\in\HeckeSym(\,\bbS(V))$ corresponds to a triple $(t,g,q)\in P$ then the 
Hecke symmetry $(\ze\ot\ze)\circ R\circ(\ze\ot\ze)^{-1}$ corresponds to 
$\ze\cdot(t,g,q)$. Therefore $R$ commutes with $\ze\ot\ze$ if and only if 
$(t,g,q)\in P(\ze)$. Such Hecke symmetries are in a bijective correspondence 
with the $\bbk^\times$-orbits in $P(\ze)$. Composing this bijection with the 
twisting transformation $R\mapsto R_\ze$ we get a bijection asserted in 
Proposition 4.1.

For any $\ph\in\GL(V)$ the space $\{\ze(x)y-\ze(y)x\mid x,y\in V\}$ of 
quadratic defining relations of the algebra $\bbS(V)_\ze$ is stable under 
the linear operator $\ph\ot\ph$ if and only if 
$\ph{\mskip1mu}\ze\ph^{-1}=c\mskip1mu\ze$ for some $c\in\bbk^\times$. In other 
words, $\ph$ extends to an automorphism of $\bbS(V)_\ze$ if and only if 
$\ph\in G(\ze)$. This means that the equivalence classes of Hecke symmetries 
in the set $\,\HeckeSym(\,\bbS(V)_\ze)$ are precisely the orbits with respect 
to the conjugating action of the group $G(\ze)$ (see the discussion at the 
end of section 1).

If $\ph\in G(\ze)$, then $(\ph\ot\ph)(\ze\ot\Id_V)=c(\ze\ot\Id_V)(\ph\ot\ph)$ 
for some $c\in\bbk^\times$, and it follows that
$$
(\ph\ot\ph)\circ R_\ze\circ(\ph\ot\ph)^{-1}=R'_\ze\quad\Lrar\quad
(\ph\ot\ph)\circ R\circ(\ph\ot\ph)^{-1}=R'.
$$
Hence the constructed bijection between 
$\,\HeckeSym(\,\bbS(V)_\ze)\setm\{{R_0}_\ze\}$ and $P(\ze)/\mskip1mu\bbk^\times$ 
is $G(\ze)$-equivariant. We thus get the second assertion of Proposition 4.1.
\endproof

Further on we assume that $\bbk$ is algebraically closed. Note that for 
$c\in\bbk^\times$ the action of the scalar operator $c\cdot\Id_V\in\GL(V)$ on 
the set $P$ coincides with the action of the element $c^2\in\bbk^\times$ 
defined by formula (4.6). It follows that each $\bbk^\times$-orbit in $P$ is 
contained in a $G(\ze)$-orbit, and therefore the $G(\ze)$-orbits in 
$P(\ze)/\mskip1mu\bbk^\times$ are in a bijection with those in $P(\ze)$.

The set $P(\ze)$ can be determined for each $\ze\in\GL(V)$, and this leads to 
the classification of the corresponding Hecke symmetries. In this paper we 
will investigate in detail the case of a diagonalizable twisting operator.

Equivalence classes in the set $\,\HeckeSym(\,\bbS(V))$ are of 8 types described 
in \cite{Sk23b}. Types 1 and 2 include Hecke symmetries with parameter $q\ne1$, 
while $q=1$ in the other types. The action of a Hecke symmetry $R$ of respective 
type is given with respect to a suitable basis $x_1,x_2,x_3$ of $V$ by the 
following formulas:
$$
\hfuzz=1.8em
\vcenter{\openup1\jot
\halign{\hfil$#$&${}#$\hfil&
\qquad\hfil$#$&${}#$\hfil&\qquad\hfil$#$&${}#$\hfil\cr
\multispan2{\hbox{\bf Type 1.}}\cr
R(x_1^2)&=qx_1^2&R(x_1x_2)&=(q-1)x_1x_2+x_2x_1&R(x_1x_3)&=(q-1)x_1x_3+x_3x_1\cr
R(x_2x_1)&=qx_1x_2&R(x_2^2)&=qx_2^2&R(x_2x_3)&=qx_3x_2\cr
R(x_3x_1)&=qx_1x_3&R(x_3x_2)&=(q-1)x_3x_2+x_2x_3&R(x_3^2)&=qx_3^2-x_1x_2+x_2x_1\cr
}}
$$
{\bf Type 2.}
The same formulas as in Type 1 with the exception that $\,R(x_3^2)=qx_3^2$.
$$
\hfuzz=1.8em
\vcenter{\openup1\jot
\halign{\hfil$#$&${}#$\hfil&
\qquad\hfil$#$&${}#$\hfil&\qquad$#$\hfil\cr
\multispan2{\hbox{\bf Type 3.}}\cr
R(x_1^2)&=x_1^2+x_1x_2-x_2x_1&R(x_1x_2)&=x_2x_1&R(x_1x_3)=x_3x_1-x_2x_3+x_3x_2\cr
R(x_2x_1)&=x_1x_2&R(x_2^2)&=x_2^2&R(x_2x_3)=x_3x_2\cr
R(x_3x_1)&=x_1x_3-x_2x_3+x_3x_2&R(x_3x_2)&=x_2x_3&R(x_3^2)=x_3^2+2(x_1x_3-x_3x_1)\cr
}}
$$
{\bf Type 4.}
As in Type 3 with the exception that $\,R(x_3^2)=x_3^2-x_1x_2+x_2x_1$.

\smallskip
\noindent
{\bf Type 5.}
As in Type 3 with the exception that $\,R(x_3^2)=x_3^2$.
$$
\hfuzz=1.8em
\vcenter{\openup1\jot
\halign{\hfil$#$&${}#$\hfil&
\qquad\hfil$#$&${}#$\hfil&\qquad$#$\hfil\cr
\hbox{\bf Type 6.}\qquad\qquad
R(x_1^2)&=x_1^2&R(x_1x_2)&=x_2x_1&R(x_1x_3)=x_3x_1\cr
R(x_2x_1)&=x_1x_2&R(x_2^2)&=x_2^2&R(x_2x_3)=x_3x_2\cr
R(x_3x_1)&=x_1x_3&R(x_3x_2)&=x_2x_3&R(x_3^2)=x_3^2+2(x_1x_3-x_3x_1)\cr
}}
$$
{\bf Type 7.}
The same formulas as in Type 1, but with $q=1$.

\smallskip
\noindent
{\bf Type 8.}
$R$ is the flip operator $R_0$ sending $x_ix_j$ to $x_jx_i$.

\medskip
Such an operator $R$ corresponds to a pair $(t,g)$ in which $t=x_1\wedge x_2$, 
while the bilinear form $g$ depends on the type to which $R$ belongs (see 
\cite{Sk23b}).

When forming the twist $R_\ze$, the condition that $R$ should commute with 
$\ze\ot\ze$ forces the basis vectors $x_1,x_2,x_3$ to be adapted to the 
twisting operator $\ze$ in some way. Some types of Hecke symmetries may not be 
permitted by a particular twisting operator.

If $\ze$ is a scalar operator, then $\bbS(V)_\ze\cong\bbS(V)$, and twisting by 
$\ze$ does not change Hecke symmetries. The case of nonscalar diagonalizable 
operators is described as follows:

\proclaim
Proposition 4.2.
Let $\ze\in\GL(V)$ be a diagonalizable linear operator with at least two 
distinct eigenvalues. Any Hecke symmetry in the set 
$\,\HeckeSym(\,\bbS(V)_\ze)$ is the $\ze$-twist of a Hecke symmetry $R$ whose 
action is given by the formulas for one of Types\/ $1$--$8$ with respect to 
some basis of\/ $V\!$ consisting of eigenvectors $x_1,x_2,x_3$ of\/ $\ze$.

Moreover, with\/ $\al_i$ being the eigenvalue of\/ $\ze$ corresponding to the 
eigenvector $x_i,$
$$
\eqalign{
&\hbox{Types $1$ and $7$ occur only when $\al_3^2=\al_1\al_2$},\cr
&\hbox{Type $3$ does not occur},\cr
&\hbox{Type $4$ occurs only when $\al_1=\al_2=-\al_3$},\cr
&\hbox{Type $5$ occurs only when $\al_1=\al_2$},\cr
&\hbox{Type $6$ occurs only when $\al_1=\al_3$}.
}
$$
\endproclaim

\Proof.
By Theorem 1.1 any Hecke symmetry in the set $\,\HeckeSym(\,\bbS(V)_\ze)$ is 
$R_\ze$ for some $R\in\,\HeckeSym(\,\bbS(V))$ such that $R$ commutes with 
$\ze\ot\ze$. If $R$ is the flip $R_0$, then its matrix form does not depend on 
the choice of a basis for $V$. So we may assume that $R\ne R_0$. Let 
$(t,g,q)\in P(\ze)$ be the corresponding triple.

There is a 2-dimensional subspace $V(t)$ of $V$ spanned by any pair of vectors 
$a,b$ such that $t=a\wedge b$. The condition that $\ze\cdot t=c{\mskip1mu}t$ 
for some $c\in\bbk^\times$ implies that $V(t)$ is invariant under $\ze$. Since 
$\ze$ is diagonalizable, there exist its eigenvectors $x_1,x_2\in V(t)$ such 
that $t=x_1\wedge x_2$. Then $c=\al_1\al_2$ where $\al_1,\al_2$ are the 
respective eigenvalues. The condition that $\ze\cdot g=c^{-1}g$ derived from 
(4.6) means that $g(u,v)=0$ whenever $u,v$ are two eigenvectors for $\ze$ with 
eigenvalues $\la,\mu$ such that $\la\mu\ne\al_1\al_2$.

Let $\al_3$ be the third eigenvalue of $\ze$ corresponding to any eigenvector 
$x_3$ linearly independent from $x_1$ and $x_2$. There are several cases 
distinguished by some conditions on the triple of eigenvalues 
$\al_1,\al_2,\al_3$.

If $\al_1,\al_2,\al_3$ are pairwise distinct, then none of the products 
$\al_1^2$, $\al_2^2$, $\al_1\al_3$, $\al_2\al_3$ is equal to $\al_1\al_2$. In 
this case the matrix of the bilinear form $g$ with respect to the basis 
$x_1,x_2,x_3$ is
$$
\pmatrix{0&\be&0\cr \be&0&0\cr 0&0&\ga}\eqno(4.8)
$$
for some $\be,\ga\in\bbk$. In view of the previously mentioned condition on 
$g$ we must have $\ga=0$ unless $\al_3^2=\al_1\al_2$. By (4.4) and (4.5) we 
have $(q-1)^2=4\be^2$. Hence $\be$ is either $(q-1)/2$ or $(1-q)/2$. Since 
$t=-x_2\wedge x_1$, we can replace $x_1,x_2$ by the pair $-x_2,\,x_1$, and 
then $\be$ in the matrix of $g$ will be changed to $-\be$. Thus we can always 
find a basis of $V$ consisting of eigenvectors of $\ze$ with respect to which 
the matrix of $g$ will have $\be=(q-1)/2$. If $\ga\ne0$, then we can achieve 
$\ga=1$, replacing $x_3$ with its scalar multiple. If $\ga=0$, then $\be\ne0$ 
since $g\ne0$. Computing the action of $R$ defined in (4.1) we 
obtain the formulas given for Types 1, 2, or 7.

Suppose now that $\al_1=\al_2\ne\al_3$. In this case $V(t)$ is a 2-dimensional 
eigenspace of $\ze$ which is orthogonal to $x_3$ with respect to $g$ because 
$\al_1\al_3\ne\al_1^2$. If the restriction of $g$ to $V(t)$ is nondegenerate, 
then there is a basis $x_1,x_2$ for $V(t)$ consisting of isotropic with 
respect to $g$ vectors. With this choice the matrix of $g$ will have the form 
(4.8), and we continue as in the previous case. So we do also when $g$ vanishes 
on $V(t)\times V(t)$. If the restriction of $g$ to $V(t)$ has rank 1, then we 
adjust the choice of $x_1$ and $x_2$ to obtain the following matrix of $g\,$:
$$
\pmatrix{1&0&0\cr0&0&0\cr0&0&\ga}.\eqno(4.9)
$$
If $\ga\ne0$, then we achieve $\ga=1$ by scaling. This is only possible when 
$\al_3^2=\al_1^2$, i.e., $\al_3=-\al_1$. The corresponding operator $R$ acts 
by formulas of Type 4 or 5 depending on whether $\ga$ is 1 or 0.

Suppose finally that $\al_1\ne\al_2$, but $\al_3$ equals either $\al_1$ 
or $\al_2$. Since $\al_i^2\ne\al_1\al_2$, we have $g(x_i,x_i)=0$ for each $i$. 
Moreover, $g$ has zero restriction to the 2-dimensional eigenspace $U$ of 
$\ze$ corresponding to the eigenvalue $\al_3$. If $g(x_1,x_2)\ne0$, then we 
can find $x_3\in U$ orthogonal to both $x_1$ and $x_2$. In such a basis the 
matrix of $g$ is of form (4.8) which has been considered already. If 
$g(x_1,x_2)=0$, then we may assume $\al_3=\al_1$, replacing $x_1,x_2$ with the 
pair $-x_2$, $x_1$ in the case when $\al_3=\al_2$. With this assumption we 
will have $g(x_3,x_2)\ne0$, and scaling the vectors brings the matrix of $g$ 
to the form
$$
\pmatrix{0&0&0\cr0&0&1\cr0&1&0}.\eqno(4.10)
$$
In this case the action of $R$ is given by the formulas of Type 6.
\endproof

The explicit formulas for the twisted operator $R_\ze$ contain the same 
monomials as the formulas for $R$, but altered coefficients. For example, 
if $R$ is of Type 1 or 2, then
$$
R_\ze(x_1x_2)=(q-1)x_1x_2+\al_2\al_1^{-1}x_2x_1.
$$
We do not write out such formulas in full because of space considerations.

It should be stressed that the eigenvectors $x_1,x_2,x_3$ in Proposition 4.2 
are not fixed for the given operator $\ze$ but each time are suited to a 
particular Hecke symmetry. This means that different Hecke symmetries in the 
description of that proposition require different numbering of eigenspaces 
and eigenvalues. For example, if $\ze$ has an eigenvalue $\al_1$ of 
multiplicity 1 and another eigenvalue $\al_2$ of multiplicity 2, then the 
formulas for Type 5 should be used with respect to an ordered triple of 
eigenvectors with respective eigenvalues $\al_2,\al_2,\al_1$, while the 
formulas for Type 6 require a reordering, so that the respective eigenvalues 
should be $\al_2,\al_1,\al_2$.

The question concerning equivalence of twisted Hecke symmetries can be solved 
with the aid of Proposition 4.1. For this we need to know the group $G(\ze)$.

\proclaim
Lemma 4.3.
The group $G(\ze)$ either coincides with the centralizer of\/ $\ze$ in $\GL(V)$ 
or contains this centralizer as a subgroup of index $3$. In the latter case\/ 
$\chr\bbk\ne3$ and the linear operator $\ze$ has $3$ distinct characteristic 
roots $\al_1,\al_2,\al_3$ such that
$$
\al_1\al_2^{-1}=\al_2\al_3^{-1}=\al_3\al_1^{-1}=\ep\eqno(4.11)
$$
where $\ep$ is a primitive cube root of $1$.
\endproclaim

\Proof.
If $\ph\in G(\ze)$, then $\ph\mskip1mu\ze\ph^{-1}=c\mskip1mu\ze$ for 
some $c\in\bbk^\times$. If $c\ne1$, then $\ze$ maps the generalized eigenspace 
of $\ze$ corresponding to some eigenvalue $\al$ to another generalized 
eigenspace which corresponds to the eigenvalue $c\mskip1mu\al$. It follows 
that $\al$ and $c\mskip1mu\al$ have the same multiplicity as roots of the 
characteristic polynomial of $\ze$. Since $\dim V=3$, there can occur at most 
one eigenvalue of multiplicity larger than 1. If $\al$ is such an eigenvalue, 
then $c\mskip1mu\al=\al$, whence $c=1$. In the other case $\ze$ has $3$ distinct 
eigenvalues which are permuted by the operator of multiplication by $c$. Hence 
$c^3=1$, and so $c$ is a primitive cube root of $1$ whenever $c\ne1$.
\endproof

Now define the type of a triple $(t,g,q)\in P$ as the type of the 
corresponding Hecke symmetry with the associated algebra $\bbS(V)$. 
Let $\ze\in\GL(V)$ be a diagonalizable linear operator. Note that its 
centralizer $C(\ze)$ in $\GL(V)$ consists of all invertible linear operators 
which leave stable each eigenspace of $\ze$.

In Proposition 4.2 we have seen that for each $(t,g,q)\in P(\ze)$ there exist 
3 linearly independent eigenvectors $x_1,x_2,x_3$ of $\ze$ such that 
$t=x_1\wedge x_2$ and the matrix of $g$ with respect to $x_1,x_2,x_3$ has one 
of several possible forms depending on $q$ and certain relations between 
the eigenvalues of $\ze$. If $q\ne1$, then the condition that the matrix of 
$g$ is (4.8) with $\be=(q-1)/2$ determines such an ordered basis $x_1,x_2,x_3$ 
of $V$ uniquely up to scaling of vectors.

Various choices of eigenvectors $x_1,x_2,x_3$ and permissible matrices of a 
symmetric bilinear form give in this way all triples in the set $P(\ze)$.

Suppose first that the three eigenvalues of $\ze$ are pairwise distinct. Then 
$P(\ze)$ has precisely six $\bbk^\times$-orbits of Type 2 with any fixed 
$q\ne1$, each corresponding to one of 6 possible orderings of the eigenvalues. 
Each of these orbits is invariant under the action of the centralizer $C(\ze)$. 
This gives 6 elements of the set $P(\ze)/\mskip1mu\bbk^\times$ fixed by the action 
of $C(\ze)$. If $\al_1,\al_2,\al_3$ satisfy (4.11), then the group $G(\ze)$ 
contains transformations which permute cyclically the eigenspaces of $\ze$. In 
this case the 6 just mentioned elements form two $G(\ze)$-orbits in 
$P(\ze)/\mskip1mu\bbk^\times$ with 3 elements in each. They correspond to two 
equivalence classes in the set $\,\HeckeSym(\,\bbS(V)_\ze)$ with 3 different 
Hecke symmetries in each. If (4.11) does not hold, then $G(\ze)=C(\ze)$, whence 
we get 6 equivalence classes with only one Hecke symmetry in each.

If there is an eigenvalue $\la$ of $\ze$ such that $\la^2$ equals the product 
of the two other eigenvalues, then the set $P(\ze)$ contains also triples of 
Type 1 and 7 with $\al_3=\la$ for the respective choice of $x_1,x_2,x_3$. 
Since the action of $C(\ze)$ allows arbitrary scaling of the value 
$g(x_3,x_3)$ while leaving $x_1$ and $x_2$ unchanged, each $C(\ze)$-orbit of 
Type 1 or 7 in $P(\ze)/\mskip1mu\bbk^\times$ has infinitely many elements.

If (4.11) does not hold, then there is only one eigenvalue $\la$ with the 
property needed. In this case the set $P(\ze)$ has precisely two 
$C(\ze)$-orbits of Type 1 with any fixed $q\ne1$. For $q=1$ there is 
only one $C(\ze)$-orbit of Type 7 since the corresponding triple $(t,g,q)$ 
does not change when the basis $x_1,x_2,x_3$ is replaced by $-x_2,x_1,x_3$. 
If the three eigenvalues of $\ze$ satisfy (4.11), then each of them can be 
taken as $\al_3$ to obtain triples of Type 1 and 7. There are six $C(\ze)$-orbits 
of Type 1 with any $q\ne1$ and three $C(\ze)$-orbits of Type 7, while extra 
transformations in the group $G(\ze)$ permute these orbits in cycles of length 
3. Thus, regardless of (4.11), the set $P(\ze)$ always has two $G(\ze)$-orbits 
of Type 1 with any $q\ne1$ and one $G(\ze)$-orbit of Type 7 with $q=1$.

Suppose now that $\ze$ has an eigenvalue $\al_1$ of multiplicity 1 and 
an eigenvalue $\al_2$ of multiplicity 2. In this case 
$G(\ze)=C(\ze)\cong\bbk^\times\times\GL_2$. In correspondence with 3 possible 
orderings of the triple of eigenvalues $\al_1,\al_2,\al_2$, the set $P(\ze)$ 
has three $G(\ze)$-orbits of Type 2 with any fixed $q\ne1$. If $\al_1=-\al_2$, 
there is one $G(\ze)$-orbit of Type 1 with any fixed $q\ne1$. It consists of 
all elements $(t,g,q)\in P(\ze)$ such that the bivector $t$ corresponds to the 
2-dimensional eigenspace of $\ze$ and the bilinear form $g$ is nondegenerate. 
Considering all triples in $P(\ze)$ corresponding to the Hecke symmetries of 
Type either 5 or 6, as discussed in the proof of Proposition 4.2, we see that 
they form one $G(\ze)$-orbit, and so do the triples of Type 4 and Type 7 in 
the case when $\al_1=-\al_2$.

Each of these orbits of Type 1, 2, 4, 5, 6, and 7 is an algebraic variety of 
dimension at least 2, and therefore it is the union of infinitely many 
$\bbk^\times$-orbits. Hence in this case each $G(\ze)$-orbit in
$P(\ze)/\mskip1mu\bbk^\times$ has infinitely many elements.

\medskip
In Corollary 4.4 we will summarize the preceding conclusions. 
Note that for a diagonalizable linear operator $\ze$ with eigenvalues 
$\al_1,\al_2,\al_3$ the algebra $\bbS(V)_\ze$ has defining relations
$$
x_3x_2=p_1x_2x_3,\qquad x_1x_3=p_2x_3x_1,\qquad x_2x_1=p_3x_1x_2\eqno(4.12)
$$
where
$$
p_1=\al_2\al_3^{-1},\qquad p_2=\al_3\al_1^{-1},\qquad p_3=\al_1\al_2^{-1}.\eqno(4.13)
$$
Conversely, the skew polynomial algebra with defining relations (4.12) is 
isomorphic to the algebra $\bbS(V)_\ze$ for a suitable $\ze$ provided that 
$\,p_1p_2p_3=1$.

Any cyclic permutation of parameters $p_1,p_2,p_3$ results in an isomorphic 
algebra. For this reason we do not mention cyclically permuted triples of 
parameters in the list of conditions given in the next corollary.

\proclaim
Corollary 4.4.
Let $A$ be the graded algebra with generators $x_1,x_2,x_3$ and defining 
relations $(4.12)$ where $p_1p_2p_3=1$. For each nonzero $q\in\bbk$ the set\/ 
$\HeckeSym(A)$ contains finitely many equivalence classes of Hecke symmetries 
with the chosen parameter $q$. The number of equivalence classes depends on $q$ 
and $p_1,p_2,p_3$ as shown in the table below:\vadjust{\vskip-2pt}
$$
\def\htt{ & height6pt & & \cr}
\offinterlineskip
\vcenter{\halign{\hbox{#}\ \hfil&\vrule#&\hfil\quad$#$\hfil&\hfil\quad$#$\hfil\cr
& & q\ne1 & q=1 \cr
\htt
\noalign{\hrule}
\htt
$p_1,p_2,p_3$ are pairwise distinct and $p_i\ne1$ for each $i$ & & 6 & 1 \cr
\htt
$p_1=p_2\ne p_3$ and $p_i\ne1$ for each $i$ & & 8 & 2 \cr
\htt
$p_1=p_2=p_3=\ep$ where $\ep$ is a primitive cube root of\/ $1$ & & 4 & 2 \cr
\htt
$p_1\ne p_2,$\quad $p_3=1$ & & 3 & 3 \cr
\htt
$p_1=p_2=-1,$\quad $p_3=1$ & & 4 & 5 \cr
\htt
$p_1=p_2=p_3=1$ & & 2 & 6 \cr
}}
$$
\smallskip\noindent
Only ${R_0}_\ze$ and the twists of Hecke symmetries of Type $2$ in the case 
when $p_i\ne\nobreak1$ for each $i$ form equivalence classes of finite 
cardinality. In particular, $\HeckeSym(A)$ contains finitely many 
Hecke symmetries with some fixed value of $q$ only when $p_1,p_2,p_3$ are 
pairwise distinct and $p_i\ne1$ for each $i$.
\endproclaim

\Proof.
If $i,j,k$ are three distinct elements of the set $\{1,2,3\}$, then it follows 
from (4.13) that $p_i=1$ if and only if $\al_j=\al_k$, and also $p_j=p_k$ 
if and only if $\al_i^2=\al_j\al_k$. By Proposition 4.1 the equivalence classes
in the set $\HeckeSym(A)$ are in a bijective correspondence with 
the $G(\ze)$-orbits in the set $P(\ze)/\mskip1mu\bbk^\times$. So we can refer 
to the counting of orbits in the discussion preceding Corollary 4.4. The first 
three lines of the table correspond to the case of pairwise distinct 
eigenvalues $\al_1,\al_2,\al_3$. The last line records the 8 types of Hecke 
symmetries with the associated algebra $\bbS(V)$.
\endproof

Motivated by Corollary 4.4 we are led to ask

\proclaim
Question 4.5.
Let $A$ be an arbitrary graded Artin-Schelter regular algebra of global 
dimension $3$ with quadratic defining relations. Is it always true that for 
each nonzero $q\in\bbk$ there are at most finitely many equivalence classes of 
Hecke symmetries with the chosen parameter $q$ and the property that 
$\,\bbS(V,R)\cong A\,?$
\endproclaim

It is very unlikely that finiteness of this kind can be always satisfied in 
dimensions larger than 3.

\references
\nextref And17
\auth{N.,Andruskiewitsch}
\paper{An introduction to Nichols algebras}
\InBook{Quantization, Geometry and Noncommutative Structures in Mathematics and Physics}
\publisher{Springer}
\Year{2017}
\Pages{135-195}

\nextref Ar-TV90
\auth{M.,Artin;J.,Tate;M.,Van den Bergh}
\paper{Some algebras associated to automorphisms of elliptic curves}
\InBook{The Grothendieck Festschrift, Volume I}
\publisher{Birkh\"auser}
\Year{1990}
\Pages{33-85}

\nextref Ar-TV91
\auth{M.,Artin;J.,Tate;M.,Van den Bergh}
\paper{Modules over regular algebras of dimension $3$}
\journal{Invent. Math.}
\Vol{106}
\Year{1991}
\Pages{335-388}

\nextref Ew-O94
\auth{H.,Ewen;O.,Ogievetsky}
\paper{Classification of the $GL(3)$ quantum matrix groups}\hfill\break
\hbox{}\hfill arXiv:9412009.

\nextref Gur90
\auth{D.I.,Gurevich}
\paper{Algebraic aspects of the quantum Yang-Baxter equation\inRus}
\journal{Algebra i Analiz}
\Vol{2} (4)
\Year{1990}
\Pages{119-148};
\etransl{Leningrad Math. J.}
\Vol{2}
\Year{1991}
\Pages{801-828}

\nextref Hua21
\auth{H.,Huang;C.,Van Nguyen;Ch.,Ure;K.B.,Vashaw;P.,Veerapen;X.,Wang}
\paper{Twisting of graded quantum groups and solutions to the quantum Yang-Baxter equation}
arXiv:2109.11585.

\nextref Hua22
\auth{H.,Huang;C.,Van Nguyen;Ch.,Ure;K.B.,Vashaw;P.,Veerapen;X.,Wang}
\paper{Twisting Manin's universal quantum groups and comodule algebras}
arXiv:2209.11621.

\nextref Kl-Sch
\auth{A.,Klimyk;K.,Schmudgen}
\book{Quantum Groups and their Representations}
\publisher{Springer}
\Year{1997}

\nextref Lyu86
\auth{V.V.,Lyubashenko}
\paper{Hopf algebras and vector symmetries\inRus}
\journal{Uspekhi Mat. Nauk}
\Vol{41} (5)
\Year{1986}
\Pages{185-186};
\etransl{Russian Math. Surveys}
\Vol{41}
\Year{1986}
\Pages{153-154}

\nextref Mo05
\auth{S.,Montgomery}
\paper{Algebra properties invariant under twisting}
\InBook{Hopf Algebras in Noncommutative Geometry and Physics}
\publisher{Dekker}
\Year{2005}
\Pages{229-243}

\nextref Shi23
\auth{N.A.,Shishmarov}
\paper{Hecke symmetries associated with Artin-Schelter regular algebras of type $E$ and $H$}
\journal{Lobachevskii J. Math.}
\Vol{44}
\Year{2023}
\Pages{4565-4579}

\nextref Sk23a
\auth{S.,Skryabin}
\paper{Hecke symmetries: an overview of Frobenius properties}
\journal{Selecta Math.}
\Vol{29}
\Year{2023}
paper no. 35.

\nextref Sk23b
\auth{S.,Skryabin}
\paper{Hecke symmetries associated with the polynomial algebra in $3$ commuting indeterminates}
\journal{J.~Algebra}
\Vol{629}
\Year{2023}
\Pages{1-20}

\nextref Zh96
\auth{J.J.,Zhang}
\paper{Twisted graded algebras and equivalences of graded categories}
\journal{Proc. London Math. Soc.}
\Vol{72}
\Year{1996}
\Pages{281-311}

\endreferences

\bigskip
{
\reffonts\rm
\noindent
Institute of Mathematics and Mechanics,
Kazan Federal University,\hfil\break
Kremlevskaya St.~18, 420008 Kazan, Russia\par
\smallskip
\noindent
Email addresses:\quad
nashishmarov@yandex.ru,\quad
serge.skryabin@kpfu.ru\par
}

\bye